\documentclass[english,russian]{amsart}
\usepackage{amssymb,amsmath}

\newtheorem{thm}{Theorem}[section]
\newtheorem{crl}[thm]{Corollary}
\newtheorem{prp}[thm]{Proposition}

\newcommand{\der}{\partial}

\newcommand{\eps}{\varepsilon}

\newcommand{\sudda}[1]{}

\newcommand{\limline}{\limits_{-\!-\!-\!-\!-\!-\!-\!-\!-\!-\!-\!-}}
\newcommand{\limsim}{\limits_{\sim\!\sim\!\sim\!\sim\!\sim\!\sim\sim\!\sim\!\sim\!\sim\!\sim\!\sim}}
\newcommand{\limeq}{\limits_{=\!=\!=\!=\!=\!=}}

\newcommand{\limcong}{\limits_{\cong\!\cong\!\cong\!\cong\!\cong\!\cong}}
\newcommand{\limapprox}{\limits_{\approx\!\approx\!\approx\!\approx\!\approx\!\approx}}
\newcommand{\limsmile}{\limits_{\smile\!\smile\!\smile\!\smile\!\smile}}
\newcommand{\lam}{\lambda}

\title{Algebras constructed by Rota-Baxter operators}
\author{A.S. Dzhumadil'daev }

\address
{Institute of Mathematica, Pushkin street 125,  Almaty,  050000,
Kazakhstan} 


\begin{document}
\maketitle

\begin{abstract} 
For associative commutative algebras  $A$ with Rota-Baxter operator  $R$
identities of the algebra $AR=(A,\circ),$ 
 where  $a\circ b= aR(b),$ are found.
\end{abstract}

\section{Introduction} 

Let $\mathcal A=Var(f_1,\ldots,f_s)$ be variety of algebras over a field $K$ generated by polynomial identities $f_1=0,\ldots, f_s=0.$ Here $f_i=f_i(t_1,\ldots, t_{k_i}), 1\le i\le s,$ are non-commutative non-asscoiative polynomials with variables $t_1,t_2,\ldots .$ Recall that a polynomial $f_i$ is called 
polynomial identity, shortly identity,  on algebra $A$, if 
 $f_i(a_1,\ldots,a_{k_i})=0$ for any $a_1,\ldots,a_{k_i}\in A,$ where calculation of 
 $ f_i(a_1,\ldots,a_{k_i})$ are made in terms of multiplication in $A.$ We will write  $f_i=0$ is identity on $A$ if  $f_i(a_1,\ldots,a_{k_i})=0$ is identity on $A.$ Variety of algebras $Var(f_1,\ldots,f_s)$ is a class of algebras generated by identities $f_1=0,\ldots, f_s=0.$

Let $R: A\rightarrow A$ be  a linear operator. Denote by $AR$ an algebra with linear space $A$ under multiplication 
$$a\circ b=a R(b).$$
An endomorphism $R:A\rightarrow A$ is called {\it  Baxter operator with weight $\lambda\in K,$ }
(\cite{baxter}, \cite{Rota}, \cite{Guo}) , with weight $\lambda,$
 if 
$$R(a)R(b)=R(a R(b)+R(a)b+\lambda a b),$$
for any $a,b\in A.$

Denote by $\mathcal{AR}_\lambda$ class of algebras $AR,$ where $A\in {\mathcal A}$ and $R\in {\mathcal R}$ is a Baxter operator with weight $\lambda.$ Here we denote by ${\mathcal R}_\lambda$ set of Baxter operators with weight $\lambda.$
We do not know whether $\mathcal{AR}_\lambda$ generates a variety of algebras. Despite of that we try to find polynomial identities that satisfy algebras of the class  $\mathcal{AR}_\lambda.$ 

Let $K$ be a field where $\sqrt{\alpha}\in K$ for any $\alpha\in K.$ Let $AR=(A,\circ)$ and $AR'=(A,\circ')$ be RB algebras, where $a\circ b=aR(b)$ and $a\circ'b=aR'(b). $ Here $R$ is  RBO with weight $\lambda \ne 0$ and  $R':A\rightarrow A$ be linear operator defined by $R'=\eps R,$
 where $\eps=\lambda'\lambda^{-1},$ $\lambda'\ne 0.$ 
 
 Let us check that a linear map  $f:A\rightarrow A$ defined by  $f(a)=\sqrt \eps\; a$ gives us an isomorphism of  algebras $AR$ and $AR'$ and that $R'$ is RBO with  weight $\lambda'.$
 
 Since $\lambda\ne 0,$ the map $f$ is bijection between algebras $AR$ and $AR'.$ Moreover, 
$$f(a)\circ f(b)=\eps a\circ b=a\circ R'(b)=a\circ 'b,$$ 
for any $a,b\in A.$ So, $f:AR\rightarrow AR'$ is isomorphism of algebras. 

Note that
$$R'(a)R'(b)=\lambda^2 R(a)R(b)=\lambda^2 R(aR(b)+R(a)b+\eps ab)=$$
$$\lambda R(a\, \lambda R(b)+\lambda R(a)\,b+\lambda\eps\, ab)=R'(aR'(b)+R'(a)b+\eps\lambda \,ab)=
R'(aR'(b)+R'(a)b+\lambda' \,ab).$$
So, $R'$ is RBO with non-zero weight $\lambda'.$ 

In particular, if the algebra $AR$ satisfies some polynomial identity, then the algebra $AR'$ satisfies the  same polynomial identity. So, the problem on finding  polynomial identities  of the class  $\mathcal{AR}_\lambda$ can be reduced to the problem  on finding polynomial identities  of the class  $\mathcal{AR}_{\lambda_0}$ for some fixed $\lambda_0\ne 0.$ Let us take $\lambda_0=-1.$

 There are two well known RBO with weight $-1.$ 
The first one is identity operator and the second one is partial sums operator (see example 1.1.5 of \cite{Guo}). Any algebra $A$ is Rota-Baxter algebra $AR$ for identity operator $R=id.$ 

We use the following strategy. We construct three polynomial identities $rcom, f_4,f_5$ of degrees $3,4,5$ that valid for any algebra $\mathcal{AR}_\lambda,$ where $\lambda\ne 0.$ Then we show that these identities are minimal and no other identities arec there if degree is no more $5.$ These facts mean that the class of algebras ${\mathcal A}sComR_\lambda,$ where $\lambda\ne 0,$ satisfies the identities $rcom=0,f_4=0, f_5=0$ and these identities are basic identities in a space of multi-linear identities of degree no more $5.$ 

In our paper we assume that ${\mathcal A}$ is a variety of associative commutative algebras $\mathcal{A}sCom.$
 For any associative commutative algebra $A$ and for any linear operator $R:A\rightarrow A$  the algebra  $AR$ satisfies  {\it right-commutative} identity $rcom=0,$ where  
 $$rcom(t_1,t_2,t_3)=(t_1t_2)t_3-(t_1t_3)t_2.$$
 Indeed,
 $$(a\circ b)\circ c-(a\circ c)\circ b=(aR(b))R(c)-(aR(c))R(b)=
a(R(b)R(c))-a(R(c)R(b))=0,$$
for any $a,b,c\in A.$

Denote by $AR^-$ and by $AR^+$ 
the algebra $AR$ under Lie commutator and Jordan commutator respectively,
$$[a,b]=a\circ b-b\circ a,$$
$$\{a,b\}=a\circ b+b\circ a.$$
Let  $$(a,b,c)=a\circ (b\circ c)-(a\circ b)\circ c$$  
be an {\it associator.} Denote by  $\langle a,b,c\rangle$ an associator under Jordan product
$$\langle a,b,c\rangle=\{a,\{b,c\}\}-\{\{a,b\},c\}.$$

If $A$ is endowed by multiplication $\star$ and $f=f(t_1,\ldots,t_k)$ is a (non-commutative non-associative) polynomial, then notation $f^\star(a_1,\ldots,a_k)$ will mean that instead of parameters $t_1,\ldots,t_k$ we make substitutions by elements $a_1,\ldots,a_k$ of $A$ and calculate 
$f(a_1,\ldots,a_k)$ in terms of multiplication $\star$. Notation
$f^{+\star}(a_1,\ldots,a_k)$ will means that the element 
$f(a_1,\ldots,a_k)$ is calculated in  terms of Jordan commutator of the multiplication $\star$. Simalarly, notation $f^{-\star}(a_1,\ldots,a_k)$
means that the element 
$f(a_1,\ldots,a_k)$ is calculated in  terms of Lie commutator of the multiplication $\star$.
If we use multiplication $\circ$, 
then we will write $f(a_1,\ldots,a_k)$ and 
$f^\pm(a_1,\ldots,a_k)$
instead of $f^\circ(a_1,\ldots,a_k)$
and  $f^{\pm\circ}(a_1,\ldots,a_k)$
For example, 
$$\langle a,b,c\rangle=(a,b,c)^+=(a,b,c)^{+\circ}.$$

In general the right-commutativity  identity is not unique identity for the algebra $AR.$
Let us give examples of additional identities for algebra $AR.$

{\bf Example.} If $R=\der$ is derivation,
$$\der(ab)=\der(a)b+a\der(b), \qquad \forall a,b\in A,$$
then the algebra $AR$ gives us so called {\it left-Novikov} algebra \cite{Novikov}: a right-commutative  algebra with left-symmetric identity for associators
$$(a,b,c)-(b,a,c)=0, \qquad \forall a,b,c\in A.$$
Any Novikov algebra  under Lie commutator is Lie. 
Under Jordan commutator we obtain so-called {\it Tortken} algebra
\cite{tortken},
$$\{\{a,b\},\{c,d\}\}- \{\{a,d\},\{c,b\}\}-\{ \langle a,b,c\rangle, d\}+
\{ \langle a,d,c\rangle, b\}=0$$

{\bf Example.} If $A=K[x]$ is a polynomials algebra with integration operator $R=\int, $ then the right-commutativity identity is not minimal for the algebra $AR.$ It is a consequence of so called {\it left-Zinbiel} identity \cite{Zinbiel}
$$(a\circ b)\circ c=a\circ (b\circ c)+a\circ(c\circ b), \qquad \forall a,b,c\in A.$$
Any Zinbiel algebra under Jordan commutator is associative commutative \cite{Zinbiel}.
Under Lie commutator it became so-called {\it Tortkara} algebra \cite{tortkara}
$$[[a,b],[c,d]]+[[a,d],[c,b]]-[jac(a,b,c),d]-[jac[a,d,c],b]=0,$$
where 
$$jac(a,b,c)=[[a,b],c]+[[b,c],a]+[[c,a],b]$$
is Jacobian. Sometimes Zinbiel algebras are called Leibniz dual algebras or chronological
algebras \cite{Agr1}, \cite{Agr2}, \cite{Kaw}. 
If weight vanishes, Rota-Baxter operator became integration operator, and $AR$ is Zinbiel, 
$AR^+$  associative-commutative and $AR^-$  Tortkara algebras. 

In our paper we construct identities of the algebra $AR,$ where $R$ is Rota-Baxter operator with non-zero weight. Call for such $R$ the algebra $AR$ as {\it RB-algebra}. 

We construct several identities for RB-algebras. Let us define  non-commutative non-associative polynomials corresponding these identities.

The first one is right-commutative polynomial mentioned above,
$$rcom(t_1,t_2,t_3)=(t_1t_2)t_3-(t_1t_3)t_2.$$
The second one is an identity of degree $4,$ defined by 
$$f_4(t_1,t_2,t_3,t_4)=t_1([t_2,t_3]t_4)-(t_1,t_2,t_3t_4)+(t_1,t_3,t_2t_4)$$
The identities $rcom(a,b,c)=0$ and $f_4(a,b,c,d)=0$ are basic multi-linear identities for RB-algebras.

The identity $f_4(a,b,c,d)=0$ implies the identities (2) and (3) (in notations of \cite{tortkara}). The polynomial corresponding the identity (2) of \cite{tortkara} we 
denote as follows 
$$f_4'(t_1,t_2,t_3,t_4)=(t_1,t_2,[t_3,t_4])+(t_1,t_3,[t_4,t_2])+(t_1,t_4,[t_2,t_3]).
$$
By Theorem 2.3 of \cite{tortkara} any algebra with identities $rcom(a,b,c)=0, f_4(a,b,c,d)=0$ under Lie commutator is Tortkara. In other words,
$$\{rcom,f_4\}\Rightarrow tortkara^{-}$$

The third polynomial has degree $5$ and is defined by 
$$f_5(t_1,t_2,t_3,t_4,t_5)=$$ $$
(t_1,t_4,(t_2,t_5,t_3))-(t_1,t_5,(t_2,t_4,t_3))+
(t_2,t_4,(t_1,t_5,t_3))-(t_2,t_5,(t_1,t_4,t_3)).$$

We show that the identity $f_4=0$ implies the identities
$tortkara^{-}=0$ and $f_5^{+}=0$.
But the identity $f_4'=0$ does not imply the identity $f_5^{+}=0.$

Counter-example gives us the algebra 
$A=K[x]$ with multiplication
$$a\star b= a\int\!\!\int b.$$
It satisfies the identity ${f_4'}^\star=0$
but $f_5^{+\star}=0$ is not identity. For example,
$$f_5^{+\star}(1,x,x^2,x^3,x^4)=
-\frac{4537 x^{18}}{6107270400}\ne 0.$$
Note that ${f_4'}^\star=0$ is not identity also. For example,
$${f_4'}^\star(1,x,x^2,x^3)=
\frac{x^{12}}{1064448}\ne 0.$$

In \cite{Guo}, Theorem 1.1.17,  by Rota-Baxter algebra one understands an algebra
with so called double product
$$a\star_R b=R(a)R(b) .$$
We do not know whether the double product algebra satisfies a polynomial identity that does not follow from commutativity identiy. We have checked that until degree $5$ any polynomial  identity of double product algebra $A=K[x]$ with 
$$a\star b=\int\!a\,\int\!b$$
is a consequence of commutativity one. 

It is interesting to consider the following deformation of the double product  algebra under deformed product 
$$a\star_{R,\eps} b=R(a)R(b)+\eps \,a R^2(b).$$
Note that if $\eps=1, \lambda=0$ the algebra $(A,\star_{R,1})$ becames Zinbiel. 

 Our main result is the following

\begin{thm}
Let $A$ be associative commutative algebra over a field of characteristic $0$ and $R:A\rightarrow A$ is Rota-Baxter operator with non-zero weight $\lambda.$ Then the algebra $AR=(A,\circ),$ where $a\circ b=a\,R(b),$ satisfies the following identities 
$$rcom(a,b,c)=(a\circ b)\circ c-(a\circ c)\circ b=0,$$ 
$$f_4(a,b,c,d)=a\circ ([b,c]\circ d)-(a,b,c\circ d)+(a,c,b\circ d)=0.$$
Any identity of the algebra $AR$ of degree no more than $4$ are consequences of the identities $rcom=0$ and $f_4=0.$

The algebra  $AR$ under Lie commutator $AR^-=(A, [\;,\;])$ satisfies the identity
$$tortkara^{-}(a,u,b,v)=
[[a,u],[b,v]]+[[a,v],[b,u]]-[jac(a,u,b),v]-[jac(a,v,b),u]=0,
$$
where
$$jac(a,b,c)=[[a,b],c]+[[b,c],a]+[[c,a],b], \quad [a,b]=a\circ b-b\circ a.$$
Any multi-linear identity of the algebra $AR^-$ of degree no more than $4$ is a consequence of identities $com=0$ and $tortkara=0.$

The algebra  $AR$ under Jordan commutator $AR^+=(A, \{\;,\,\})$ satisfies the identity
$$f_5^{+}(a,b,c,u,v)=\langle a,u,\langle b,v,c\rangle\rangle-\langle a,v,\langle b,u,c\rangle\rangle+\langle b,u,\langle a,v,c\rangle\rangle-\langle b,v,\langle a,u,c\rangle\rangle=0.$$
Any multi-linear identity of the algebra $AR^-$ of degree no more than $5$ is a consequence of identities $acom=0$ and $f_5=0.$
 \end{thm}

Additional results concern Rota Baxter operators with weight $0,$ i.e., for integration operators on 
$A=K[x].$  We show that the multiplication defined by 
\begin{itemize}
\item $a\star_{0,n} b=\sum_{i=0}^n {n\choose i}\int_i a \int_{n-i} b$ is associative
\item $a\star_{1,n}b=\sum_{i=0}^n {n\choose i} \int_i a \int_{n-i+1}b$ is left-zinbiel
\item $[a,b]_{n}=\sum_{i=0}^n {n\choose i}\left(1-\frac{2i}{n}\right)\int_i a\int_{n-i}b$ is Tortkara
\item $a\star_{k,n} b=\sum_{i=0}^n {n\choose i} \int_i a \int_{n-i+k}b$ is right-commutative for $k>2.$
\end{itemize}
Here we define $\int $ by $\int a(x)=\int_0^x a(x)dx$ and 
$\int_i a=\int\cdots \int a$
(number of integrals is $i$)

\section{Identities of RB-algebras of degree 3 and 4}
A linear map  $R : A\rightarrow A$ is called a Baxter operator with weight $\lambda$  if it satisfies the following identity
\begin{equation}\label{Baxter}
R(a)R(b) = R(a R(b) + b R(a) + \lambda  a b),
\end{equation}
for any $a,b\in A. $

In this section we establish the following facts.

Let $A$ be Baxter algebra with Baxter operator $R.$ Then the algebra $AR = (A,\circ)$ satisfies the following identities
\begin{equation}\label{Baxter1}
(x\circ a)\circ b=(x\circ b)\circ a,
\end{equation}
\begin{equation}\label{Baxter2}
a\circ ([x,y]\circ b)-(a,x,y\circ b)+(a,y,x\circ b)=0,
\end{equation}
for any $a,b,x,y,\in A,$ where $[a,b]=a\circ b-b\circ a.$

{\bf Proof.} Relation (\ref{Baxter1}) is evident
$$(x\circ a)\circ b=(x\,R(a))\,R(b)=
(x\,R(b))\,R(a)=
(x\circ b)\circ a.$$

Let us prove (\ref{Baxter2}). We have
$$a\circ ([x,y]\circ b)-(a,x,y\circ b)+(a,y,x\circ b)=$$

$$a\circ ((x\circ y)\circ b)-a\circ ((y\circ x)\circ b)-$$ 
$$a\circ (x\circ (y\circ b))+(a\circ x)\circ (y\circ b)+$$
$$a\circ (y\circ (x\circ b))-
(a\circ y)\circ (x\circ b)=$$

$$a\,R((x R(y))R(b))-a\,R((y R(x))R(b))-
a\,R(x\,R(y\,R(b)))+(a\,R(x))\, R(y\,R(b))$$
$$+a\,R(y\,R(x\,R(b)))-(a\,R(y))\,R(x\,R(b))=
$$
$$a(Q_1+Q_2+Q_3),$$
where
$$Q_1= R(x)\, R(y\,R(b))-R((y R(x))R(b)),$$
$$Q_2=R((x R(y))R(b))-R(y)\,R(x\,R(b)),$$
$$Q_3=-R(x\,R((y\,R( b)))+R(y\,R(x\,R(b))).$$

By (\ref{Baxter}) we have 
$$Q_1=$$
$$R(x)\, R(y\,R(b))-R(R(x) (y\,R(b)))=$$
$$R(x\,R(y\,R(b)))+\lambda\,R(x\,y\,R(b)),$$

$$Q_2=$$
$$ R((x\,R(b))\,R(y))-R(y)\,R(x\,R(b))= 
$$
$$-R(y\,R(x\,R(b)))-\lambda\,R(y\,x\,R(b)).$$
Thus,
$$Q_1+Q_2=R(x\,R(y\,R(b)))-R(y\,R(x\,R(b)))=-Q_3,$$
and,
$$a\circ ([x,y]\circ b)-(a,x,y\circ b)+(a,y,x\circ b)=Q_1+Q_2+Q_3=0.$$
Identity (\ref{Baxter2}) is established. 
$\square$

\begin{crl}
If $A$ is a Baxter algebra, then the algebra
$AR = (A, \circ)$ satisfies the identity
\begin{equation}\label{Baxter3}
f_4'(x,a,b,c)= (x,a,[b,c]) + (x,b,[c,a]) + (x,c,[a,b]) = 0
\end{equation}
for any $a,b,c,x \in A.$
\end{crl}

{\bf Proof.} By  (\ref{Baxter1}) we have 
$$(x,a,[b,c])=$$ $$
x\circ (a\circ [b,c])-(x\circ a)\circ [b,c]=$$
$$
x\circ (a\circ [b,c])-(x\circ [b,c])\circ a.
$$
Similarly,
$$(x,b,[c,a])=
x\circ (b\circ [c,a])-(x\circ [c,a])\circ b,
$$

$$(x,c,[a,b])=
x\circ (c\circ [a,b])-(x\circ [a,b])\circ c.
$$
Hence
$$ (x,a,[b,c]) + (x,b,[c,a]) + (x,c,[a,b]) = $$

$$ x\circ (a\circ [b,c])-(x\circ [b,c])\circ a+
x\circ (b\circ [c,a])-(x\circ [c,a])\circ b+
x\circ (c\circ [a,b])-(x\circ [a,b])\circ c=
$$

$$-(x\circ [a,b])\circ c-(x\circ [b,c])\circ a
-(x\circ [c,a])\circ b+$$
$$ x\circ (a\circ [b,c])+
x\circ (b\circ [c,a])+
x\circ (c\circ [a,b])=$$

$$(x,[a,b],c)-x\circ ([a,b]\circ c)+$$
$$(x,[b,c],a)-x\circ ([b,c]\circ a)+$$
$$(x,[c,a],b)-x\circ ([c,a]\circ b)+$$
$$ x\circ (a\circ (b\circ c))-x\circ (a\circ (c\circ b))+$$
$$
x\circ (b\circ (c\circ a))-
x\circ (b\circ (a\circ c))+$$ $$
x\circ (c\circ (a\circ b))-x\circ (c\circ (b\circ a))=
$$

$$-x\circ ([a,b]\circ c)-x\circ ([b,c]\circ a)-x\circ ([c,a]\circ b)+$$
$$(x,[a,b],c)+(x,[b,c],a)+(x,[c,a],b])+$$
$$ x\circ (a\circ (b\circ c))-x\circ (a\circ (c\circ b))+$$
$$
x\circ (b\circ (c\circ a))-
x\circ (b\circ (a\circ c))+$$ $$
x\circ (c\circ (a\circ b))-x\circ (c\circ (b\circ a))=
$$

$$-x\circ ([a,b]\circ c)
+ x\circ (a\circ (b\circ c))
-x\circ (b\circ (a\circ c))$$
$$-x\circ ([b,c]\circ a)+x\circ (b\circ (c\circ a))-x\circ (c\circ (b\circ a))$$
$$-x\circ ([c,a]\circ b)+x\circ (c\circ (a\circ b))-x\circ (a\circ (c\circ b))+S,$$
where
$$S=(x,[a,b],c)+(x,[b,c],a)+(x,[c,a],b]),$$
We have 
$$S=$$

$$(x,a\circ b,c)-(x,b\circ a,c)
+(x,b\circ c,a)-(x,c\circ b,a)
+(x,c\circ a,b)-(x,a\circ c,b)=$$

$$x\circ((a\circ b)\circ c)
-(x\circ (a\circ b))\circ c$$
$$
-x\circ ((b\circ a)\circ c)
+(x\circ (b\circ a))\circ c$$
$$+x\circ ((b\circ c)\circ a)-
(x\circ (b\circ c))\circ a$$
$$-x\circ ((c\circ b)\circ a)
+(x\circ (c\circ b))\circ a$$
$$+x\circ ((c\circ a)\circ b)-
(x\circ (c\circ a))\circ b$$
$$-x\circ ((a\circ c)\circ b)+
+(x\circ (a\circ c))\circ b=$$

$$x\circ((a\circ b)\circ c)
-x\circ ((a\circ c)\circ b)$$
$$-x\circ ((b\circ a)\circ c)
+x\circ ((b\circ c)\circ a)$$
$$-x\circ ((c\circ b)\circ a)
+x\circ ((c\circ a)\circ b)$$
$$
-(x\circ (a\circ b))\circ c
+(x\circ (b\circ a))\circ c
-(x\circ (b\circ c))\circ a
+(x\circ (c\circ b))\circ a
-(x\circ (c\circ a))\circ b
+(x\circ (a\circ c))\circ b=$$
(by (\ref{Baxter1}))
$$
-(x\circ c)\circ (a\circ b)
+(x\circ c)\circ (b\circ a)
-(x\circ a)\circ (b\circ c)
+(x\circ a)\circ (c\circ b)
-(x\circ b)\circ(c\circ a)
+(x\circ b)\circ (a\circ c)$$

Therefore,
$$ (x,a,[b,c]) + (x,b,[c,a]) + (x,c,[a,b]) = $$
$$-x\circ ([a,b]\circ c)
+ x\circ (a\circ (b\circ c))
-x\circ (b\circ (a\circ c))$$
$$-x\circ ([b,c]\circ a)+x\circ (b\circ (c\circ a))-x\circ (c\circ (b\circ a))$$
$$-x\circ ([c,a]\circ b)+x\circ (c\circ (a\circ b))-x\circ (a\circ (c\circ b))$$
$$
-(x\circ c)\circ (a\circ b)
+(x\circ c)\circ (b\circ a)
-(x\circ a)\circ (b\circ c)
+(x\circ a)\circ (c\circ b)
-(x\circ b)\circ(c\circ a)
+(x\circ b)\circ (a\circ c)=$$

$$-x\circ ([a,b]\circ c)
+ x\circ (a\circ (b\circ c))-(x\circ a)\circ (b\circ c)
-x\circ (b\circ (a\circ c))
+(x\circ b)\circ (a\circ c)
$$
$$-x\circ ([b,c]\circ a)+x\circ (b\circ (c\circ a))-(x\circ b)\circ(c\circ a)-x\circ (c\circ (b\circ a))
+(x\circ c)\circ (b\circ a)
$$
$$-x\circ ([c,a]\circ b)+x\circ (c\circ (a\circ b))-(x\circ c)\circ (a\circ b)-x\circ (a\circ (c\circ b))+
(x\circ a)\circ (c\circ b)
=$$

$$-x\circ ([a,b]\circ c)+(x,a,b\circ c)-(x,b,a\circ c)$$
$$-x\circ ([b,c]\circ a)+(x,b,c\circ a)-
(x,c,b\circ a)$$
$$-x\circ ([c,a]\circ b)+(x,c,a\circ b)-
(x,a,c\circ b).$$
Hence by (\ref{Baxter2})
$$ (x,a,[b,c]) + (x,b,[c,a]) + (x,c,[a,b]) = 0$$
is identity for any $a,b,c,c,x\in A.$
$\square$

\section{Identities of RB-algebras of degree 5}


In this section we establish the following facts.

Let $A$ be associative commutative algebra with Baxter operator $R: A\rightarrow A.$ Let us define polynomial $f(a,b,c,u,v)$ on $A$ by  
$$f(a,b,c,u,v)=
-\langle a,v,\langle b,u,c\rangle\rangle+
\langle a,u,\langle b,v,c\rangle\rangle,
$$
where
$$\langle a,b,c\rangle=\{a,\{b,c\}\}-\{\{a,b\},c\},$$
$$\{a,b\}=a\circ b+b\circ a,$$
$$a\circ b=a\,R(b).$$
Then 
$$f(a,b,c,u,v)=-f (b, a, c, u, v),$$
$$f(a,b,c,u,v)=-f(a,b,c,v,u),$$
for any $a,b,c,u,v\in A.$

{\bf Proof.} Let us define $R$-associator by 
$$(a,b,c)_R=a\,R(bc)-R(ab)\,c.$$
Because of commutativity of $A$ we see that
\begin{equation}\label{assR}
(a,b,c)_R=-(c,b,a)_R.
\end{equation}


Let us prove that
$$\langle a,b,c\rangle=-\lambda (a,b,c)_R$$
We have
$$\langle a,b,c\rangle=$$
$$\{a,\{b,c\}\}-\{\{a,b\},c\}=$$
$$a\circ (b\circ c)+a\circ (c\circ b)+
(b\circ c)\circ a+(c\circ b)\circ a$$
$$-(a\circ b)\circ c-(b\circ a)\circ c-c\circ (a\circ b)-c\circ (b\circ a)=$$
(by (\ref{Baxter1}) )

$$a\circ (b\circ c)+a\circ (c\circ b)+(c\circ b)\circ a$$
$$-(a\circ b)\circ c-c\circ (a\circ b)-c\circ (b\circ a)=$$

$$a\,R(b\,R(c))+a\,R (c\,R(b))+(c\,R( b)\,R( a)
-(a\,R( b))\,R( c)-c\,R(a\,R( b))-c\,R (b\,R( a))=$$
(by (\ref{Baxter}))

$$a\,R(b)\,R(c))-\lambda a\,R(bc) +c\,R(a)\, R(b)
-a\,R( b)\,R( c)-c\,R(a)R( b)+\lambda c\,R(ab)
=$$

$$\lambda(- a\,R(bc) + c\,R(ab))=$$
$$-\lambda(a,b,c)_R.$$
Therefore,
$$f(a,b,c,u,v)=$$

$$\langle a,v,\lambda(b\,R(uc)-R(b u)\,c)\rangle -
\langle a,u,\lambda(b\,R(vc)-R(v b)\,c)\rangle=
$$

$$\lambda \langle  a,v,b\,R(uc)\rangle -
\lambda \langle a,v, 
R(b u)\,c)\rangle -\lambda \langle  a,u,b\,R(vc)\rangle +
\lambda \langle a,u, 
R(b v)\,c)\rangle=$$

$$-\lambda^2 a\,R(v\,b\,R(uc))+\lambda^2 R(av) b\,R(uc)+\lambda^2 a\,R(c v R(b u))
-\lambda^2 R(a v) R(b u)c $$
$$
+\lambda^2 a\,R(u\,b\,R(vc))-\lambda^2 R(au) b\,R(vc)-\lambda^2 a\,R(c u R(b v))
+\lambda^2 R(a u) R(b v)c.
$$
Hence,
$$\lambda^{-2}(f(a,b,c,u,v)+f(b,a,c,u,v))=$$

$$-a\,R(v\,b\,R(uc))+ R(av) b\,R(uc)+ a\,R(c v R(b u))- R(a v) R(b u)c $$
$$+a\,R(u\,b\,R(vc))- R(au) b\,R(vc)- a\,R(c u R(b v))+ R(a u) R(b v)c$$

$$-b\,R(v\,a\,R(uc))+ R(bv) a\,R(uc)+ b\,R(c v R(a u))- R(b v) R(au)c $$
$$+b\,R(u\,a\,R(vc))- R(bu) a\,R(vc)- b\,R(c u R(a v))+ R(b u) R(a v)c=$$

$$-a\,R(v\,b\,R(uc))+ R(av) b\,R(uc)+ a\,R(c v R(b u)) $$
$$+a\,R(u\,b\,R(vc))- R(au) b\,R(vc)- a\,R(c u R(b v))$$

$$-b\,R(v\,a\,R(uc))+ R(bv) a\,R(uc)+ b\,R(c v R(a u)) $$
$$+b\,R(u\,a\,R(vc))- R(bu) a\,R(vc)- b\,R(c u R(a v))=$$

$$-a\,R(v\,b\,R(uc))+ a\,R(c v R(b u)) $$
$$+a\,R(u\,b\,R(vc))- a\,R(c u R(b v))$$
$$+ R(bv) a\,R(uc)- R(bu) a\,R(vc)$$

$$-b\,R(v\,a\,R(uc))+ b\,R(c v R(a u)) $$
$$+b\,R(u\,a\,R(vc))- b\,R(c u R(a v))$$
$$+ R(av) b\,R(uc)- R(au) b\,R(vc)=$$

$$a[ R(bv)R(uc)-R(v\,b\,R(uc))]$$
$$+a[- R(bu)R(vc)+ R(c v R(b u))]$$
$$+a[R(u\,b\,R(vc)) - R(c u R(b v))]$$

$$+b[ R(av) R(uc)-R(v\,a\,R(uc))]$$
$$+b[- R(au) R(vc)+ R(c v R(a u))]$$
$$+b[R(u\,a\,R(vc)) - R(c u R(a v))]=$$
(by (\ref{Baxter}) )

$$a[\mathop{R(R(vb)u c)}\limline+\lambda \mathop{R(bcuv)}\limeq]-$$
$$a[\mathop{R(R(c v) b u)}\limsim+\lambda \mathop{R(bcuv)}\limeq]+$$
$$a[\mathop{R(u\,b\,R(vc))}\limsim - \mathop{R(c u R(b v))}\limline]+$$

$$b[\mathop{R(R(v\,a)uc)}\limapprox+
\lambda \mathop{R(acuv)}\limcong]-$$
$$b[\mathop{R(R(c v) a u)}\limsmile+\lambda \mathop{R(acuv)}\limcong]+$$
$$b[\mathop{R(u\,a\,R(vc))}\limsmile - \mathop{R(c u R(a v))}\limapprox]=$$

$$0.$$
$\square$

\section{Multi-linear identities of $AR^-$ in degree $4$\label{four}}

Let us write skew-symmetric non-associative  polynomial in degree $4$ 
in terms of commutators.

There are two bracketing  types of skew-symmetric monomials for four variables: type $((ab)c)d$ and type $(ab)(cd).$ In terms of non-planar binary  trees these two types can be  constructed by the following trees 

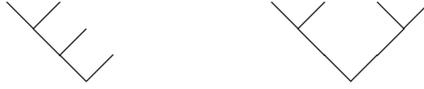
\begin{picture}(100,100)

\put(130,30){\line(-1,1){30}}
\put(130,30){\line(1,1){10}}
\put(120,40){\line(1,1){10}}
\put(110,50){\line(1,1){10}}

\put(230,30){\line(-1,1){30}}
\put(230,30){\line(1,1){30}}
\put(210,50){\line(1,1){10}}
\put(250,50){\line(-1,1){10}}
\end{picture}

So, general form of such polynomial in degree $4$ is 
$$X_{4,com}^-(a,b,c,d)=$$
$$\lambda_{1} [[[a, b], c], d] + 
 \lambda_{2} [[[a, b], d], c] + \lambda_{3} [[[a, c], b], d] + 
 \lambda_{4} [[[a, c], d], b] + \lambda_{5} [[[a, d], b], c] + $$ $$
 \lambda_{6} [[[a, d], c], b] + \lambda_{7} [[[b, c], a], d] + 
 \lambda_{8} [[[b, c], d], a] + \lambda_{9} [[[b, d], a], c] +  
 \lambda_{10} [[[b, d], c], a] + $$ $$\lambda_{11} [[[c, d], a], b] + 
 \lambda_{12} [[[c, d], b], a]+ \lambda_{13} [[a, b], [c, d]] + \lambda_{14} [[a, c], [b, d]] + 
 \lambda_{15} [[a, d], [b, c]]. $$

\bigskip

The aim of this section is to show that any identity of algebra $AR^-$ in  degree $4$ follows from the identity $tortkara^{lie}(a,b,c,d)=0$

To do that we use realisation of Rota-Baxter operator $R$ as operator on sequences given by (see \cite{Guo})
$$R((a_1,a_2,a_3,a_4,\ldots ))=(a_1,a_1+a_2,a_1+a_2+a_3,a_1+a_2+a_3+a_4,\ldots).$$
Multiplication of sequences is component wise. 

For us is enough to collect  relations  that give fourth components of operator $R.$ Below we use  notations of the form 
$a=(0,1,0,1)$ instead of $a=(0,1,0,1,0,0,\ldots).$

Let  $X_{4,com}^-(a,b,c,d)$ be an  element of $AR^-$ calculated in terms of Lie commutator.
For example, if 
$$a=(0, 1, 0, 1), b=(1, 0, 1, 0), c=(0, 1, 1, 1), d=(1, 1, 1, 1),$$
then first two components of $X_{4,com}^-(a,b,c,d)$ vanishes,
$$(X_{4,com}^-(a,b,c,d))_i=0,\quad i=1,2,$$
and third component is 
$$ (X_{4,com}^-(a,b,c,d))_3=$$
$$
-4 \lambda_{1} - 4 \lambda_{2} - 2 \lambda_{3} - 2 \lambda_{4} - 3 \lambda_{5} - 
  2 \lambda_{6} + \lambda_{8} + \lambda_{9} + 2 \lambda_{10} + \lambda_{11} - 2 \lambda_{13} + 
  \lambda_{14} + \lambda_{15}.$$
 But we pay attention  for the fourth component only 
 $$ (X_{4,com}^-(a,b,c,d))_4=$$ 
$$
 14 \lambda_{1} + 14 \lambda_{2} + 7 \lambda_{3} + 8 \lambda_{4} + 9 \lambda_{5} + 8 \lambda_{6} - 
  3 \lambda_{7} - 5 \lambda_{8} - 5 \lambda_{9} - 6 \lambda_{10} $$ $$- 2 \lambda_{11} + \lambda_{12} + 
  4 \lambda_{13} - 2 \lambda_{14} - 2 \lambda_{15}.$$
In the Table 1 we give twelve relations that give fourth components of $X_{4,com}^-(a,b,c,d)$ where $a=(0,1,0,1)$ and sequences $b,c,d$ are given below

\bigskip

{\bf Table 1.}

$$
\begin{array}{|c|c|c|l|}
\hline
b&c&d&\mbox{fourth component of $X_{4,com}^-(a,b,c,d)$}\\
\hline\hline
  (1, 0, 1, 0)& (0, 1, 1, 1)& (1, 1, 1, 1)& 
  14 \lambda_{1} + 14 \lambda_{2} + 7 \lambda_{3} + 8 \lambda_{4} + 9 \lambda_{5} +
   8 \lambda_{6} - 
   3 \lambda_{7} - 5 \lambda_{8} - 5 \lambda_{9} -\\&&& 6 \lambda_{10} -2 \lambda_{11} + \lambda_{12} + 
   4 \lambda_{13} - 2 \lambda_{14} - 
   2 \lambda_{15}\\
   \hline
 (1, 0, 1, 0)& (0, 1, 1, 1)& (1, 1, 2, 1)& 
  18 \lambda_{1} + 19 \lambda_{2} + 9 \lambda_{3} + 10 \lambda_{4} + 14 \lambda_{5} + 
   14 \lambda_{6} - 4 \lambda_{7} - 8 \lambda_{8} - 2 \lambda_{9} - \\&&&4 \lambda_{10} + 
   2 \lambda_{11} + 4 \lambda_{12} + 2 \lambda_{13} - 3 \lambda_{14} - 
   5 \lambda_{15}\\
   \hline

  (1, 0, 1, 0)& (0, 1, 1, 1)& (1, 2, 1, 1)& 
  18 \lambda_{1} + 19 \lambda_{2} + 9 \lambda_{3} + 10 \lambda_{4} + 13 \lambda_{5} + 
   12 \lambda_{6} - 4 \lambda_{7} - 7 \lambda_{8} - 6 \lambda_{9} -\\
   &&& 8 \lambda_{10} - 
   2 \lambda_{11} + 2 \lambda_{12} + 6 \lambda_{13} - 2 \lambda_{14} - 
   3 \lambda_{15}\\
   \hline
   
  (1, 0, 1, 0)& (0, 1, 1, 2)& (1, 1, 2, 1)& 
  18 \lambda_{1} + 22 \lambda_{2} + \lambda_{3} + 2 \lambda_{4} + 16 \lambda_{5} + 16 \lambda_{6} - 
   12 \lambda_{7} - 16 \lambda_{8} - 2 \lambda_{9} - \\&&&3 \lambda_{10} + 10 \lambda_{11} + 
   12 \lambda_{12} + 2 \lambda_{13} - 2 \lambda_{14} - 
   7 \lambda_{15}\\
   \hline
  (1, 0, 1, 0)& (0, 1, 1, 2)& (1, 2, 1, 1)& 
  18 \lambda_{1} + 22 \lambda_{2} + \lambda_{3} + 2 \lambda_{4} + 14 \lambda_{5} + 12 \lambda_{6} - 
   12 \lambda_{7} - 15 \lambda_{8} - 8 \lambda_{9} - \\&&&8 \lambda_{10} + 6 \lambda_{11} + 
   10 \lambda_{12} + 6 \lambda_{13} - 2 \lambda_{14} - 
   3 \lambda_{15}\\
   \hline
  (1, 0, 1, 0)& (0, 1, 2, 1)& (1, 1, 1, 1)& 
  21 \lambda_{1} + 20 \lambda_{2} + 14 \lambda_{3} + 16 \lambda_{4} + 13 \lambda_{5} + 
   12 \lambda_{6} + \lambda_{7} - 2 \lambda_{8} - 7 \lambda_{9} -\\&&& 9 \lambda_{10} - 8 \lambda_{11} - 
   3 \lambda_{12} + 8 \lambda_{13} - 4 \lambda_{14} - 
   4 \lambda_{15}\\ 
   \hline
   (1, 0, 1, 0)& (0, 1, 2, 1)& (1, 2, 1, 1)& 
  27 \lambda_{1} + 27 \lambda_{2} + 18 \lambda_{3} + 20 \lambda_{4} + 19 \lambda_{5} + 
   18 \lambda_{6} + \lambda_{7} - 3 \lambda_{8} - 8 \lambda_{9} -\\&&& 12 \lambda_{10} - 
   10 \lambda_{11} - 3 \lambda_{12} + 12 \lambda_{13} - 4 \lambda_{14} - 
   6 \lambda_{15}\\
   \hline
   (1, 0, 1, 0)& (0, 1, 2, 2)& (1, 1, 2, 1)& 
  27 \lambda_{1} + 30 \lambda_{2} + 10 \lambda_{3} + 12 \lambda_{4} + 22 \lambda_{5} + 
   22 \lambda_{6} - 7 \lambda_{7} - 13 \lambda_{8} - 3 \lambda_{9} - \\&&&6 \lambda_{10} + 
   4 \lambda_{11} +8 \lambda_{12} + 6 \lambda_{13} - 5 \lambda_{14} - 
   10 \lambda_{15}\\ \hline 
   (1, 0, 1, 0)& (0, 1, 2, 2)& (1, 2, 2, 1)& 
  33 \lambda_{1} + 38 \lambda_{2} + 12 \lambda_{3} + 14 \lambda_{4} + 28 \lambda_{5} + 
   28 \lambda_{6} - 9 \lambda_{7} - 16 \lambda_{8} - 5 \lambda_{9} -\\&&& 9 \lambda_{10} + 
   4 \lambda_{11} + 10 \lambda_{12} + 10 \lambda_{13} - 5 \lambda_{14} - 
   12 \lambda_{15}\\ 
   \hline
   (1, 0, 2, 0)& (0, 1, 1, 1)& (1, 1, 1, 1)& 
  24 \lambda_{1} + 24 \lambda_{2} + 10 \lambda_{3} + 12 \lambda_{4} + 14 \lambda_{5} + 
   12 \lambda_{6} - 10 \lambda_{7} - 11 \lambda_{8} - 13 \lambda_{9} - \\&&& 12 \lambda_{10} -
   3 \lambda_{11} + 3 \lambda_{12} + 7 \lambda_{13} - 
   \lambda_{14}\\ \hline
   (1, 0, 2, 0)& (0, 1, 1, 1)& (1, 1, 2, 1)& 
  31 \lambda_{1} + 31 \lambda_{2} + 13 \lambda_{3} + 15 \lambda_{4} + 21 \lambda_{5} + 
   21 \lambda_{6} - 13 \lambda_{7} - 15 \lambda_{8} -\\&&& 10 \lambda_{9} - 10 \lambda_{10} + 
   3 \lambda_{11} + 7 \lambda_{12} + 5 \lambda_{13} - 2 \lambda_{14} - 
   3 \lambda_{15}\\ \hline
   (1, 0, 2, 0)& (0, 1, 1, 1)& (1, 2, 1, 1)& 
  31 \lambda_{1} + 32 \lambda_{2} + 13 \lambda_{3} + 15 \lambda_{4} + 20 \lambda_{5} + 
   18 \lambda_{6} - 13 \lambda_{7} - 15 \lambda_{8} -\\&&& 17 \lambda_{9} - 16 \lambda_{10} - 
   3 \lambda_{11} + 5 \lambda_{12} + 11 \lambda_{13}\\
   \hline
\end{array}$$

These relations gives us linear system of $12$ equations and $15$ unknowns $\lambda_i, 1\le i\le 15.$
This system has rank $12$ and unknowns $\lambda_1, \lambda_2, \lambda_4$ can be selected as a free parameters. 
Solution of this system:
$$\lambda_{3} = -\lambda_{1}, $$ $$\lambda_{5} = -\lambda_{2}, $$ $$\lambda_{6} = -\lambda_{4}, $$ $$
  \lambda_{7} = \lambda_{1},$$ $$ \lambda_{8} = \lambda_{1} - \lambda_{2} + \lambda_{4}, $$ $$
  \lambda_{9} = \lambda_{2}, $$ $$ \lambda_{10} = -\lambda_{1} + \lambda_{2} - \lambda_{4}, $$ $$
  \lambda_{11} = \lambda_{4}, $$ $$ \lambda_{12} = \lambda_{1} - \lambda_{2} + \lambda_{4}, $$ $$
  \lambda_{13} = -\lambda_{1} + \lambda_{2}, $$ $$ \lambda_{14} = \lambda_{1} + \lambda_{4}, $$ $$
  \lambda_{15} = \lambda_{2} - \lambda_{4}$$
Hence, if $X_4^{com}(a,b,c,d)=0$ is identity for 
algebra $AR^{-},$ where $R$ is Baxter operator, then  
$$X_4^{com}(a,b,c,d)=$$
$$
  \sum_{i=1}^3 \lambda_{i}g_i(a,b,c,d),$$ 
  where
  $$g_1(a,b,c,d)=
  -(a b) (c d) + (a c) (b d) + 
    ((a b) c) d - ((a c) b) d +$$ $$ ((b c) a) d + 
    ((b c) d) a -
     ((b d) c) a + ((c d) b) a,$$ 
     
     $$g_2(a,b,c,d)=
(a b) (c d) + (a d) (b c) + 
    ((a b) d) c - ((a d) b) c - $$ $$((b c) d) a+ 
    ((b d) a) c + ((b d) c) a - ((c d) b) a,$$ 
    
    $$
g_3(a,b,c,d)=(a c) (b d) - (a d) (b c) + 
    ((a c) d) b - ((a d) c) b + $$ $$((b c) d) a - 
    ((b d) c) a +((c d) a) b + ((c d) b) a]=
$$

We see that
$$g_1(a,b,c,d)=tortkara(c, a, d, b),$$

$$g_2(a,b,c,d)=
tortkara(a, c, b, d) - tortkara(b, a, d, c),$$

$$g_3(a,b,c,d)=tortkara(b, a, d, c).$$
Hence any multi-linear  identity of degree $4$ of the algebra $AR^{-},$ where $R$ is Baxter operator, is a consequence of the identity $tortkara=0.$ $\square$

\section{ Multilinear right-commutative identity in degree 4 }

General form of the right-commutative polynomial of degree $4$
$$X_{4,rcom}(a, b, c, d) = $$

 $$ \lambda_{1} ((a b) c) d + \lambda_{2} ((b a) c) d + 
  \lambda_{3} ((c a) b) d + \lambda_{4} ((d a) b) c+
  \lambda_{5} a ((b c) d) + $$ $$\lambda_{6} a ((c b) d) + 
  \lambda_{7} a ((d b) c))+ \lambda_{8} b ((a c) d) + 
  \lambda_{9} b ((c a) d) + \lambda_{10} b ((d a) c))+ $$ $$
  \lambda_{11} c ((a b) d) + 
  \lambda_{12} c ((b a) d) + 
  \lambda_{13} c ((d a) b)+ \lambda_{14} d ((a b) c) + 
\lambda_{15} d ((b a) c) + $$ $$\lambda_{16} d ((c a) b) + 
 \lambda_{17} (a b) (c d) + \lambda_{18} (a b) (d c) + 
  \lambda_{19} (a c) (b d) + \lambda_{20} (a c) (d b) + $$ $$
  \lambda_{21} (a d) (b c) +\lambda_{22} (a d) (c b) + 
  \lambda_{23} (b a) (c d) + \lambda_{24} (b a) (d c) + 
  \lambda_{25} (b c) (a d) + $$ $$\lambda_{26} (b c) (d a) + 
  \lambda_{27} (b d) (a c) + \lambda_{28} (b d) (c a) + 
  \lambda_{29} (c a) (b d) + \lambda_{30} (c a) (d b) + $$ $$
  \lambda_{31} (c b) (a d) + 
  \lambda_{32} (c b) (d a) + 
  \lambda_{33} (c d) (a b) + \lambda_{34} (c d) (b a) + 
  \lambda_{35} (d a) (b c) + $$ $$\lambda_{36} (d a) (c b) + 
  \lambda_{37} (d b) (a c) + \lambda_{38} (d b) (c a) + 
  \lambda_{39} (d c) (a b) + \lambda_{40} (d c) (b a) +$$ $$
\lambda_{41} a (b (c d)) + 
\lambda_{42} a (b (d c)) + 
\lambda_{41} a (b (c d)) + \lambda_{42} a (b (d c)) + 
  \lambda_{43} a (c (b d)) +$$ $$ \lambda_{44} a (c (d b)) + 
  \lambda_{45} a (d (b c)) + \lambda_{46} a (d (c b)) + 
 \lambda_{47} b (a (c d)) +  \lambda_{48} b (a (d c)) + $$ $$
\lambda_{49} b (c (a d)) + 
  \lambda_{50} b (c (d a)) + \lambda_{51} b (d (a c)) + 
  \lambda_{52} b (d (c a)) +
   \lambda_{53} c (a (b d)) +$$ $$ \lambda_{54} c (a (d b)) + 
  \lambda_{55} c (b (a d)) + \lambda_{56} c (b (d a)) + 
 \lambda_{57} c (d (a b)) + \lambda_{58} c (d (b a)) + $$ $$\lambda_{59} d (a (b c)) + 
  \lambda_{60} d (a (c b)) +
 \lambda_{61} d (b (a c)) + 
  \lambda_{62} d (b (c a)) + \lambda_{63} d (c (a b)) + $$ $$
  \lambda_{64} d (c (b a)) .$$ 

One checks that $X_{4,rcom}=0$ is identity for sequences RB-algebra, if 
$$X_{4,rcom}(a,b,c,d)=\sum_{i=1}^{12} \lam_i g_i(a,b,c,d),$$
where

$$g_1(a,b,c,d)=$$ 
$$-a (b (d c)) + a (d (b c)) + a ((b c) d) - 
   a ((d b) c) + (a b) (d c) - (a d) (b c),$$

  $$g_2(a,b,c,d)=$$ $$
  -a (c (d b)) + a (d (c b)) + a ((c b) d) - 
   a ((d b) c) + (a c) (d b) - (a d) (c b),$$

   $$g_3(a,b,c,d)=$$ $$ 
  -a (b (c d)) + a (b (d c)) + a (c (b d)) - 
   a (c (d b)) - a (d (b c)) + a (d (c b)) + $$
   $$
   (a b) (c d) - (a b) (d c) - (a c) (b d) + 
   (a c) (d b) + (a d) (b c) - (a d) (c b),$$

    $$g_4(a,b,c,d)=$$ $$
  -b (a (d c)) + b (d (a c)) + b ((a c) d) - 
   b ((d a) c) + (b a) (d c) - (b d) (a c),$$

    $$g_5(a,b,c,d)=$$ $$
  -b (c (d a)) + b (d (c a)) + b ((c a) d) - 
   b ((d a) c) + (b c) (d a) - (b d) (c a),$$

    $$g_6(a,b,c,d)=$$ $$
  -b (a (c d)) + b (a (d c)) + b (c (a d)) - 
   b (c (d a)) - b (d (a c)) + b (d (c a)) + $$
   $$
   (b a) (c d) - (b a) (d c) - (b c) (a d) + 
   (b c) (d a) + (b d) (a c) - (b d) (c a),$$

    $$g_7(a,b,c,d)=$$ $$
  -c (a (d b)) + c (d (a b)) + c ((a b) d) - 
   c ((d a) b) + (c a) (d b) - (c d) (a b),$$

    $$g_8(a,b,c,d)=$$ $$
  -c (b (d a)) + c (d (b a)) + c ((b a) d) - 
   c ((d a) b) + (c b) (d a) - (c d) (b a),$$

    $$g_9(a,b,c,d)=$$ $$
  -c (a (b d)) + c (a (d b)) + c (b (a d)) - 
   c (b (d a)) - c (d (a b)) + c (d (b a)) + $$
   $$
   (c a) (b d) - (c a) (d b) - (c b) (a d) + 
   (c b) (d a) + (c d) (a b) - (c d) (b a),$$

    $$g_{10}(a,b,c,d)=$$ $$
  -d (a (c b)) + d (c (a b)) + d ((a b) c) - 
   d ((c a) b) + (d a) (c b) - (d c) (a b),$$

    $$g_{11}(a,b,c,d)=$$ $$
  -d (b (c a)) + d (c (b a)) + d ((b a) c) - 
   (d ((c a) b) + (d b) (c a) - (d c) (b a),$$

    $$g_{12}(a,b,c,d)=$$ $$
  -d (a (b c)) + d (a (c b)) + d (b (a c)) - 
   d (b (c a)) - d (c (a b)) + d (c (b a)) + $$
   $$
   (d a) (b c) - (d a) (c b) - (d b) (a c) + 
   (d b) (c a) + (d c) (a b) - (d c) (b a).$$
By right-commutative rule we see that 

 $$g_1(a,b,c,d)=f_4(a, b, d, c),$$
 $$g_2(a,b,c,d)= f_4(a, c, d, b),$$
  $$g_3(a,b,c,d)=
 f_4(a, b, c, d) - f_4(a, b, d, c) +  f_4(a, c, d, b),$$
  $$g_4(a,b,c,d)=
 f_4(b, a, d, c),$$
  $$g_5(a,b,c,d)=
 f_4(b, c, d, a),$$
  $$g_6(a,b,c,d)=
 f_4(b, a, c, d) - f_4(b, a, d, c) + f_4(b, c, d, a), $$
  $$g_7(a,b,c,d)=
f_4(c, a, d, b),$$
  $$g_8(a,b,c,d)=f_4(c, b, d, a),$$
$$g_9(a,b,c,d)= f_4(c, a, b, d) - f_4(c, a, d, b) +  f_4(c, b, d, a), $$
 $$g_{10}(a,b,c,d)=f_4(d, a, c, b),$$
  $$g_{11}(a,b,c,d)=f_4(d, b, c, a),$$
  $$g_{12}(a,b,c,d)=f_4(d, a, b, c) - f_4(d, a, c, b) + f_4(d, b, c, a).$$
Hence any identity of the algebra $AR,$ where $R$ is Rota-Baxter agebra with non-zero weight, in degree $4$ follows form the identity $f_4(a,b,c,d)=0.$
$\square$

\section{Multi-linear identities of $AR^+$ in degree $5$}

Let 
$$X_{5}(a,b,c,d,e)=$$ 
$$X_{5,1}(a,b,c,d,e)+X_{5,2}(a,b,c,d,e)+X_{5,3}(a,b,c,d,e)$$
be general form of commutative polynomial in degree $5,$
where 
$$X_{5,1}(a,b,c,d,e)=$$
$$\lambda_{1} \{\{\{\{a, b\}, c\}, d\}, e\} + \lambda_{2} \{\{\{\{a, b\}, c\}, e\}, d\} +
  \lambda_{3} \{\{\{\{a, b\}, d\}, c\}, e\} + 
 \lambda_{4} \{\{\{\{a, b\}, d\}, e\}, c\} +$$ $$ \lambda_{5} \{\{\{\{a, b\}, e\}, c\}, d\} +
  \lambda_{6} \{\{\{\{a, b\}, e\}, d\}, c\} + 
 \lambda_{7} \{\{\{\{a, c\}, b\}, d\}, e\} + \lambda_{8} \{\{\{\{a, c\}, b\}, e\}, d\} +$$ $$
  \lambda_{9} \{\{\{\{a, c\}, d\}, b\}, e\} + 
 \lambda_{10} \{\{\{\{a, c\}, d\}, e\}, b\} + 
 \lambda_{11} \{\{\{\{a, c\}, e\}, b\}, d\} + 
 \lambda_{12} \{\{\{\{a, c\}, e\}, d\}, b\} + $$ $$
 \lambda_{13} \{\{\{\{a, d\}, b\}, c\}, e\} + 
 \lambda_{14} \{\{\{\{a, d\}, b\}, e\}, c\} + 
 \lambda_{15} \{\{\{\{a, d\}, c\}, b\}, e\} +
 \lambda_{16} \{\{\{\{a, d\}, c\}, e\}, b\} + $$ $$
 \lambda_{17} \{\{\{\{a, d\}, e\}, b\}, c\} + 
 \lambda_{18} \{\{\{\{a, d\}, e\}, c\}, b\} + 
 \lambda_{19} \{\{\{\{a, e\}, b\}, c\}, d\} + 
 \lambda_{20} \{\{\{\{a, e\}, b\}, d\}, c\} + $$ $$
 \lambda_{21} \{\{\{\{a, e\}, c\}, b\}, d\} + 
 \lambda_{22} \{\{\{\{a, e\}, c\}, d\}, b\} + 
 \lambda_{23} \{\{\{\{a, e\}, d\}, b\}, c\} + 
 \lambda_{24} \{\{\{\{a, e\}, d\}, c\}, b\} + $$ $$
 \lambda_{25} \{\{\{\{b, c\}, a\}, d\}, e\} + 
 \lambda_{26} \{\{\{\{b, c\}, a\}, e\}, d\} + 
 \lambda_{27} \{\{\{\{b, c\}, d\}, a\}, e\} + 
 \lambda_{28} \{\{\{\{b, c\}, d\}, e\}, a\} + $$ $$
 \lambda_{29} \{\{\{\{b, c\}, e\}, a\}, d\} + 
 \lambda_{30} \{\{\{\{b, c\}, e\}, d\}, a\} + 
 \lambda_{31} \{\{\{\{b, d\}, a\}, c\}, e\} + 
 \lambda_{32} \{\{\{\{b, d\}, a\}, e\}, c\} + $$ $$
 \lambda_{33} \{\{\{\{b, d\}, c\}, a\}, e\} + 
 \lambda_{34} \{\{\{\{b, d\}, c\}, e\}, a\} + 
 \lambda_{35} \{\{\{\{b, d\}, e\}, a\}, c\} + 
 \lambda_{36} \{\{\{\{b, d\}, e\}, c\}, a\} + $$ $$
 \lambda_{37} \{\{\{\{b, e\}, a\}, c\}, d\} + 
 \lambda_{38} \{\{\{\{b, e\}, a\}, d\}, c\} + 
 \lambda_{39} \{\{\{\{b, e\}, c\}, a\}, d\} + 
 \lambda_{40} \{\{\{\{b, e\}, c\}, d\}, a\} + $$ $$
 \lambda_{41} \{\{\{\{b, e\}, d\}, a\}, c\} + 
 \lambda_{42} \{\{\{\{b, e\}, d\}, c\}, a\} + 
 \lambda_{43} \{\{\{\{c, d\}, a\}, b\}, e\} + 
 \lambda_{44} \{\{\{\{c, d\}, a\}, e\}, b\} + $$ $$
 \lambda_{45} \{\{\{\{c, d\}, b\}, a\}, e\} + 
 \lambda_{46} \{\{\{\{c, d\}, b\}, e\}, a\} + 
 \lambda_{47} \{\{\{\{c, d\}, e\}, a\}, b\} + 
 \lambda_{48} \{\{\{\{c, d\}, e\}, b\}, a\} + $$ $$
 \lambda_{49} \{\{\{\{c, e\}, a\}, b\}, d\} + 
 \lambda_{50} \{\{\{\{c, e\}, a\}, d\}, b\} + 
 \lambda_{51} \{\{\{\{c, e\}, b\}, a\}, d\} + 
 \lambda_{52} \{\{\{\{c, e\}, b\}, d\}, a\} + $$ $$
 \lambda_{53} \{\{\{\{c, e\}, d\}, a\}, b\} + 
 \lambda_{54} \{\{\{\{c, e\}, d\}, b\}, a\} + 
 \lambda_{55} \{\{\{\{d, e\}, a\}, b\}, c\} + 
 \lambda_{56} \{\{\{\{d, e\}, a\}, c\}, b\} + $$ $$
 \lambda_{57} \{\{\{\{d, e\}, b\}, a\}, c\} + 
 \lambda_{58} \{\{\{\{d, e\}, b\}, c\}, a\} + 
 \lambda_{59} \{\{\{\{d, e\}, c\}, a\}, b\} + \lambda_{60} \{\{\{\{d, e\}, c\}, b\}, a\},$$

\bigskip

$$X_{5,2}(a,b,c,d,e)=$$
$$
\lambda_{61} \{\{\{a, b\}, \{c, d\}\}, e\} + 
 \lambda_{62} \{\{\{a, b\}, \{c, e\}\}, d\} + 
 \lambda_{63} \{\{\{a, b\}, \{d, e\}\}, c\} + 
 \lambda_{64} \{\{\{a, c\}, \{b, d\}\}, e\} + $$ $$
 \lambda_{65} \{\{\{a, c\}, \{b, e\}\}, d\} + 
 \lambda_{66} \{\{\{a, c\}, \{d, e\}\}, b\} + 
\lambda_{67} \{\{\{a, d\}, \{b, c\}\}, e\} + 
 \lambda_{68} \{\{\{a, d\}, \{b, e\}\}, c\} + $$ $$
 \lambda_{69} \{\{\{a, d\}, \{c, e\}\}, b\} + 
\lambda_{70} \{\{\{a, e\}, \{b, c\}\}, d\} + 
 \lambda_{71} \{\{\{a, e\}, \{b, d\}\}, c\} + 
 \lambda_{72} \{\{\{a, e\}, \{c, d\}\}, b\} + $$ $$
 \lambda_{73} \{\{\{b, c\}, \{d, e\}\}, a\}+
\lambda_{74} \{\{\{b, d\}, \{c, e\}\}, a\} + 
\lambda_{75} \{\{\{b, e\}, \{c, d\}\}, a\} ,
$$

$$X_{5,3}(a,b,c,d,e)=$$
$$
\lambda_{76} \{\{\{a, b\}, c\}, \{d, e\}\} + 
 \lambda_{77} \{\{\{a, b\}, d\}, \{c, e\}\} + 
 \lambda_{78} \{\{\{a, b\}, e\}, \{c, d\}\} + 
\lambda_{79} \{\{\{a, c\}, b\}, \{d, e\}\} + $$ $$
 \lambda_{80} \{\{\{a, c\}, d\}, \{b, e\}\} + 
 \lambda_{81} \{\{\{a, c\}, e\}, \{b, d\}\} + 
\lambda_{82} \{\{\{a, d\}, b\}, \{c, e\}\} + 
 \lambda_{83} \{\{\{a, d\}, c\}, \{b, e\}\} + $$ 
 $$
 \lambda_{84} \{\{\{a, d\}, e\}, \{b, c\}\} + 
\lambda_{85} \{\{\{a, e\}, b\}, \{c, d\}\} + 
 \lambda_{86} \{\{\{a, e\}, c\}, \{b, d\}\} + 
 \lambda_{87} \{\{\{a, e\}, d\}, \{b, c\}\} + $$ 
 $$
\lambda_{88} \{\{\{b, c\}, a\}, \{d, e\}\} + 
 \lambda_{89} \{\{\{b, c\}, d\}, \{a, e\}\} + 
 \lambda_{90} \{\{\{b, c\}, e\}, \{a, d\}\} + 
 \lambda_{91} \{\{\{b, d\}, a\}, \{c, e\}\} + $$ $$
 \lambda_{92} \{\{\{b, d\}, c\}, \{a, e\}\} + 
 \lambda_{93} \{\{\{b, d\}, e\}, \{a, c\}\} + 
\lambda_{94} \{\{\{b, e\}, a\}, \{c, d\}\} + 
 \lambda_{95} \{\{\{b, e\}, c\}, \{a, d\}\} + $$ $$
 \lambda_{96} \{\{\{b, e\}, d\}, \{a, c\}\} + 
\lambda_{97} \{\{\{c, d\}, a\}, \{b, e\}\} + 
 \lambda_{98} \{\{\{c, d\}, b\}, \{a, e\}\} + 
 \lambda_{99} \{\{\{c, d\}, e\}, \{a, b\}\} + $$ $$
 \lambda_{100} \{\{\{c, e\}, a\}, \{b, d\}\} + 
 \lambda_{101} \{\{\{c, e\}, b\}, \{a, d\}\} + 
 \lambda_{102} \{\{\{c, e\}, d\}, \{a, b\}\} + 
 \lambda_{103} \{\{\{d, e\}, a\}, \{b, c\}\} + $$ $$
 \lambda_{104} \{\{\{d, e\}, b\}, \{a, c\}\} + 
 \lambda_{105} \{\{\{d, e\}, c\}, \{a, b\}\}.$$

We have 
$$f_5^+(a,b,c,d,e)=$$
$$  \langle(a,d,\langle b,e,c\rangle\rangle-\langle(a,e,\langle b,d,c\rangle\rangle
+\langle(b,d,\langle a,e,c\rangle\rangle-\langle(b,e,\langle a,d,c\rangle\rangle
.$$

Calculations similar to calculations given in section \ref{four} show that 
$$X_ {5,com}^+ (a,b,c,d,e)=\sum_{i=1}^{20}\lambda_i g_i(a,b,c,d,e),$$
where Jordan polynomials $g_i(a,b,c,d,e)$
are defined as follows: 

$$g_1(a,b,c,d,e)=$$
$$
-\{\{\{a, c\}, d\}, \{b, e\}\} + \{\{\{a, d\}, c\}, \{b, e\}\} + 
    \{\{\{b, c\}, e\}, \{a, d\}\} - \{\{\{b, d\}, e\}, \{a, c\}\} - $$ $$
    \{\{\{b, e\}, c\}, \{a, d\}\} + \{\{\{b, e\}, d\}, \{a, c\}\} + 
    \{\{\{\{a, c\}, d\}, b\}, e\} - \{\{\{\{a, d\}, c\}, b\}, e\} - $$ $$
    \{\{\{\{b, c\}, d\}, a\}, e\} + \{\{\{\{b, c\}, d\}, e\}, a\} - 
    \{\{\{\{b, c\}, e\}, d\}, a\} + \{\{\{\{b, d\}, c\}, a\}, e\} - $$ $$
    \{\{\{\{b, d\}, c\}, e\}, a\} + \{\{\{\{b, d\}, e\}, c\}, a\} + 
    \{\{\{\{b, e\}, c\}, d\}, a\} - 
    \{\{\{\{b, e\}, d\}, c\}, a\},
$$

$$g_2(a,b,c,d,e)=$$ $$
 -\{\{\{a, d\}, e\}, \{b, c\}\} + 
    \{\{\{a, e\}, d\}, \{b, c\}\} - \{\{\{b, c\}, d\}, \{a, e\}\} + 
    \{\{\{b, c\}, e\}, \{a, d\}\} +$$ $$ \{\{\{b, d\}, c\}, \{a, e\}\} - 
    \{\{\{b, e\}, c\}, \{a, d\}\} + \{\{\{\{a, d\}, e\}, b\}, c\} - 
    \{\{\{\{a, e\}, d\}, b\}, c\} +$$ $$ \{\{\{\{b, c\}, d\}, e\}, a\} - 
    \{\{\{\{b, c\}, e\}, d\}, a\} - \{\{\{\{b, d\}, c\}, e\}, a\} - 
    \{\{\{\{b, d\}, e\}, a\}, c\} + $$ $$\{\{\{\{b, d\}, e\}, c\}, a\} + 
    \{\{\{\{b, e\}, c\}, d\}, a\} + \{\{\{\{b, e\}, d\}, a\}, c\} - 
    \{\{\{\{b, e\}, d\}, c\}, a\},
$$

$$g_3(a,b,c,d,e)=$$
$$
 -\{\{\{a, c\}, e\}, \{b, d\}\} + 
    \{\{\{a, e\}, c\}, \{b, d\}\} + \{\{\{b, c\}, d\}, \{a, e\}\} - 
    \{\{\{b, d\}, c\}, \{a, e\}\} +$$ $$ \{\{\{b, d\}, e\}, \{a, c\}\} - 
    \{\{\{b, e\}, d\}, \{a, c\}\} + \{\{\{\{a, c\}, e\}, b\}, d\} - 
    \{\{\{\{a, e\}, c\}, b\}, d\} - $$ 
    $$\{\{\{\{b, c\}, d\}, e\}, a\} - 
    \{\{\{\{b, c\}, e\}, a\}, d\} + \{\{\{\{b, c\}, e\}, d\}, a\} + 
    \{\{\{\{b, d\}, c\}, e\}, a\} -$$ $$ \{\{\{\{b, d\}, e\}, c\}, a\} + 
    \{\{\{\{b, e\}, c\}, a\}, d\} - \{\{\{\{b, e\}, c\}, d\}, a\} + 
    \{\{\{\{b, e\}, d\}, c\}, a\}, 
$$

$$g_4(a,b,c,d,e)=$$
$$
-\{\{\{a, d\}, b\}, \{c, e\}\} + 
    \{\{\{a, d\}, c\}, \{b, e\}\} + \{\{\{b, d\}, a\}, \{c, e\}\} - 
    \{\{\{b, e\}, c\}, \{a, d\}\} -$$ $$ \{\{\{c, d\}, a\}, \{b, e\}\} + 
    \{\{\{c, e\}, b\}, \{a, d\}\} + \{\{\{\{a, d\}, b\}, c\}, e\} - 
    \{\{\{\{a, d\}, c\}, b\}, e\} -$$ $$ \{\{\{\{b, d\}, a\}, c\}, e\} + 
    \{\{\{\{b, d\}, c\}, a\}, e\} -\{\{\{\{b, d\}, c\}, e\}, a\} + 
    \{\{\{\{b, e\}, c\}, d\}, a\} +$$ $$ \{\{\{\{c, d\}, a\}, b\}, e\} - 
    \{\{\{\{c, d\}, b\}, a\}, e\} + \{\{\{\{c, d\}, b\}, e\}, a\} - 
    \{\{\{\{c, e\}, b\}, d\}, a\},
$$

$$g_5(a,b,c,d,e)=$$
$$
 -\{\{\{a, e\}, b\}, \{c, d\}\} + 
    \{\{\{a, e\}, c\}, \{b, d\}\} - \{\{\{b, d\}, c\}, \{a, e\}\} + 
    \{\{\{b, e\}, a\}, \{c, d\}\} +$$ $$ \{\{\{c, d\}, b\}, \{a, e\}\} - 
    \{\{\{c, e\}, a\}, \{b, d\}\} + \{\{\{\{a, e\}, b\}, c\}, d\} - 
    \{\{\{\{a, e\}, c\}, b\}, d\} +$$ $$ \{\{\{\{b, d\}, c\}, e\}, a\} - 
    \{\{\{\{b, e\}, a\}, c\}, d\} + \{\{\{\{b, e\}, c\}, a\}, d\} - 
    \{\{\{\{b, e\}, c\}, d\}, a\} -$$ $$ \{\{\{\{c, d\}, b\}, e\}, a\} + 
    \{\{\{\{c, e\}, a\}, b\}, d\} - \{\{\{\{c, e\}, b\}, a\}, d\} + 
    \{\{\{\{c, e\}, b\}, d\}, a\},
$$

$$g_6(a,b,c,d,e)=$$
$$
 -\{\{\{a, b\}, d\}, \{c, e\}\} + 
    \{\{\{a, d\}, c\}, \{b, e\}\} + \{\{\{b, c\}, e\}, \{a, d\}\} + 
    \{\{\{b, d\}, a\}, \{c, e\}\} - $$ $$\{\{\{b, e\}, c\}, \{a, d\}\} - 
    \{\{\{c, d\}, a\}, \{b, e\}\} - \{\{\{c, d\}, e\}, \{a, b\}\} + 
    \{\{\{c, e\}, d\}, \{a, b\}\} +$$ $$ \{\{\{\{a, b\}, d\}, c\}, e\} - 
    \{\{\{\{a, d\}, c\}, b\}, e\} - \{\{\{\{b, c\}, d\}, a\}, e\} + 
    \{\{\{\{b, c\}, d\}, e\}, a\} - $$ $$\{\{\{\{b, c\}, e\}, d\}, a\} - 
    \{\{\{\{b, d\}, a\}, c\}, e\} + \{\{\{\{b, d\}, c\}, a\}, e\} - 
    \{\{\{\{b, d\}, c\}, e\}, a\} +$$ $$ \{\{\{\{b, e\}, c\}, d\}, a\} + 
    \{\{\{\{c, d\}, a\}, b\}, e\} + \{\{\{\{c, d\}, e\}, b\}, a\} - 
    \{\{\{\{c, e\}, d\}, b\}, a\},
$$

$$g_7(a,b,c,d,e)=$$
$$
 -\{\{\{a, d\}, e\}, \{b, c\}\} + 
    \{\{\{a, e\}, d\}, \{b, c\}\} - \{\{\{b, c\}, d\}, \{a, e\}\} + 
    \{\{\{b, c\}, e\}, \{a, d\}\} +$$ $$ \{\{\{c, d\}, b\}, \{a, e\}\} - 
    \{\{\{c, e\}, b\}, \{a, d\}\} + \{\{\{\{a, d\}, e\}, c\}, b\} - 
    \{\{\{\{a, e\}, d\}, c\}, b\} + $$ $$\{\{\{\{b, c\}, d\}, e\}, a\} - 
    \{\{\{\{b, c\}, e\}, d\}, a\} - \{\{\{\{c, d\}, b\}, e\}, a\} - 
    \{\{\{\{c, d\}, e\}, a\}, b\} +$$ $$ \{\{\{\{c, d\}, e\}, b\}, a\} + 
    \{\{\{\{c, e\}, b\}, d\}, a\} + \{\{\{\{c, e\}, d\}, a\}, b\} - 
    \{\{\{\{c, e\}, d\}, b\}, a\},
$$

$$g_8(a,b,c,d,e)=$$
$$
 -\{\{\{a, b\}, e\}, \{c, d\}\} + 
    \{\{\{a, e\}, c\}, \{b, d\}\} + \{\{\{b, c\}, d\}, \{a, e\}\} - 
    \{\{\{b, d\}, c\}, \{a, e\}\} +$$ $$ \{\{\{b, e\}, a\}, \{c, d\}\} + 
    \{\{\{c, d\}, e\}, \{a, b\}\} - \{\{\{c, e\}, a\}, \{b, d\}\} - 
    \{\{\{c, e\}, d\}, \{a, b\}\} +$$ $$ \{\{\{\{a, b\}, e\}, c\}, d\} - 
    \{\{\{\{a, e\}, c\}, b\}, d\} - \{\{\{\{b, c\}, d\}, e\}, a\} - 
    \{\{\{\{b, c\}, e\}, a\}, d\} +$$ $$ \{\{\{\{b, c\}, e\}, d\}, a\} + 
    \{\{\{\{b, d\}, c\}, e\}, a\} - \{\{\{\{b, e\}, a\}, c\}, d\} + 
    \{\{\{\{b, e\}, c\}, a\}, d\} -$$ $$ \{\{\{\{b, e\}, c\}, d\}, a\} - 
    \{\{\{\{c, d\}, e\}, b\}, a\} + \{\{\{\{c, e\}, a\}, b\}, d\} + 
    \{\{\{\{c, e\}, d\}, b\}, a\}, 
$$

$$g_9(a,b,c,d,e)=$$
$$
-\{\{\{a, c\}, b\}, \{d, e\}\} + 
    \{\{\{a, d\}, c\}, \{b, e\}\} + \{\{\{b, c\}, a\}, \{d, e\}\} + 
    \{\{\{b, c\}, e\}, \{a, d\}\} -$$ $$ \{\{\{b, d\}, e\}, \{a, c\}\} - 
    \{\{\{b, e\}, c\}, \{a, d\}\} - \{\{\{c, d\}, a\}, \{b, e\}\} + 
    \{\{\{d, e\}, b\}, \{a, c\}\} +$$ $$ \{\{\{\{a, c\}, b\}, d\}, e\} - 
    \{\{\{\{a, d\}, c\}, b\}, e\} - \{\{\{\{b, c\}, a\}, d\}, e\} - 
    \{\{\{\{b, c\}, e\}, d\}, a\} + $$ $$\{\{\{\{b, d\}, c\}, a\}, e\} - 
    \{\{\{\{b, d\}, c\}, e\}, a\} + \{\{\{\{b, d\}, e\}, c\}, a\} + 
    \{\{\{\{b, e\}, c\}, d\}, a\} +$$ $$ \{\{\{\{c, d\}, a\}, b\}, e\} - 
    \{\{\{\{c, d\}, b\}, a\}, e\} + \{\{\{\{c, d\}, b\}, e\}, a\} - 
    \{\{\{\{d, e\}, b\}, c\}, a\},
$$

$$g_{10}(a,b,c,d,e)=$$
$$
 -\{\{\{a, c\}, b\}, \{d, e\}\} + 
    \{\{\{a, e\}, c\}, \{b, d\}\} + \{\{\{b, c\}, a\}, \{d, e\}\} + 
    \{\{\{b, c\}, d\}, \{a, e\}\} -$$ $$ \{\{\{b, d\}, c\}, \{a, e\}\} - 
    \{\{\{b, e\}, d\}, \{a, c\}\} - \{\{\{c, e\}, a\}, \{b, d\}\} + 
    \{\{\{d, e\}, b\}, \{a, c\}\} + $$ $$\{\{\{\{a, c\}, b\}, e\}, d\} - 
    \{\{\{\{a, e\}, c\}, b\}, d\} - \{\{\{\{b, c\}, a\}, e\}, d\} - 
    \{\{\{\{b, c\}, d\}, e\}, a\} +$$ $$ \{\{\{\{b, d\}, c\}, e\}, a\} + 
    \{\{\{\{b, e\}, c\}, a\}, d\} - \{\{\{\{b, e\}, c\}, d\}, a\} + 
    \{\{\{\{b, e\}, d\}, c\}, a\} +$$ $$ \{\{\{\{c, e\}, a\}, b\}, d\} - 
    \{\{\{\{c, e\}, b\}, a\}, d\} + \{\{\{\{c, e\}, b\}, d\}, a\} - 
    \{\{\{\{d, e\}, b\}, c\}, a\},
$$

$$g_{11}(a,b,c,d,e)=$$
$$
 -\{\{\{a, c\}, d\}, \{b, e\}\} + 
    \{\{\{a, e\}, d\}, \{b, c\}\} - \{\{\{b, c\}, d\}, \{a, e\}\} + 
    \{\{\{b, e\}, d\}, \{a, c\}\} +$$ $$ \{\{\{c, d\}, a\}, \{b, e\}\} + 
    \{\{\{c, d\}, b\}, \{a, e\}\} - \{\{\{d, e\}, a\}, \{b, c\}\} - 
    \{\{\{d, e\}, b\}, \{a, c\}\} +$$ $$ \{\{\{\{a, c\}, d\}, e\}, b\} - 
    \{\{\{\{a, e\}, d\}, c\}, b\} + \{\{\{\{b, c\}, d\}, e\}, a\} - 
    \{\{\{\{b, e\}, d\}, c\}, a\} - $$ $$\{\{\{\{c, d\}, a\}, e\}, b\} - 
    \{\{\{\{c, d\}, b\}, e\}, a\} + \{\{\{\{d, e\}, a\}, c\}, b\} + 
    \{\{\{\{d, e\}, b\}, c\}, a\}, 
$$

$$g_{12}(a,b,c,d,e)=$$
$$
-\{\{\{a, c\}, e\}, \{b, d\}\} + 
    \{\{\{a, e\}, d\}, \{b, c\}\} - \{\{\{b, c\}, d\}, \{a, e\}\} + 
    \{\{\{b, d\}, e\}, \{a, c\}\} +$$ $$ \{\{\{c, d\}, b\}, \{a, e\}\} + 
    \{\{\{c, e\}, a\}, \{b, d\}\} - \{\{\{d, e\}, a\}, \{b, c\}\} - 
    \{\{\{d, e\}, b\}, \{a, c\}\} +$$ $$ \{\{\{\{a, c\}, e\}, d\}, b\} - 
    \{\{\{\{a, e\}, d\}, c\}, b\} + \{\{\{\{b, c\}, d\}, e\}, a\} - 
    \{\{\{\{b, d\}, e\}, c\}, a\} -$$ $$ \{\{\{\{c, d\}, b\}, e\}, a\} - 
    \{\{\{\{c, d\}, e\}, a\}, b\} + \{\{\{\{c, d\}, e\}, b\}, a\} - 
    \{\{\{\{c, e\}, a\}, d\}, b\} + $$ $$\{\{\{\{c, e\}, d\}, a\}, b\} - 
    \{\{\{\{c, e\}, d\}, b\}, a\} + \{\{\{\{d, e\}, a\}, c\}, b\} + 
    \{\{\{\{d, e\}, b\}, c\}, a\},
$$

$$g_{13}(a,b,c,d,e)=$$
$$
 -\{\{\{a, e\}, b\}, \{c, d\}\} + 
    \{\{\{a, e\}, d\}, \{b, c\}\} - \{\{\{b, c\}, d\}, \{a, e\}\} + 
    \{\{\{b, e\}, a\}, \{c, d\}\} +$$ $$ \{\{\{c, d\}, b\}, \{a, e\}\} - 
    \{\{\{d, e\}, a\}, \{b, c\}\} + \{\{\{\{a, e\}, b\}, d\}, c\} - 
    \{\{\{\{a, e\}, d\}, b\}, c\} +$$ $$ \{\{\{\{b, c\}, d\}, e\}, a\} - 
    \{\{\{\{b, e\}, a\}, d\}, c\} + \{\{\{\{b, e\}, d\}, a\}, c\} - 
    \{\{\{\{b, e\}, d\}, c\}, a\} -$$ $$ \{\{\{\{c, d\}, b\}, e\}, a\} + 
    \{\{\{\{d, e\}, a\}, b\}, c\} - \{\{\{\{d, e\}, b\}, a\}, c\} + 
    \{\{\{\{d, e\}, b\}, c\}, a\},
$$

$$g_{14}(a,b,c,d,e)=$$
$$
 -\{\{\{a, d\}, b\}, \{c, e\}\} + 
    \{\{\{a, e\}, d\}, \{b, c\}\} - \{\{\{b, c\}, d\}, \{a, e\}\} + 
    \{\{\{b, d\}, a\}, \{c, e\}\} +$$ $$ \{\{\{b, d\}, c\}, \{a, e\}\} - 
    \{\{\{b, e\}, c\}, \{a, d\}\} + \{\{\{c, e\}, b\}, \{a, d\}\} - 
    \{\{\{d, e\}, a\}, \{b, c\}\} +$$ $$ \{\{\{\{a, d\}, b\}, e\}, c\} - 
    \{\{\{\{a, e\}, d\}, b\}, c\} + \{\{\{\{b, c\}, d\}, e\}, a\} - 
    \{\{\{\{b, d\}, a\}, e\}, c\} -$$ $$ \{\{\{\{b, d\}, c\}, e\}, a\} + 
    \{\{\{\{b, e\}, c\}, d\}, a\} + \{\{\{\{b, e\}, d\}, a\}, c\} - 
    \{\{\{\{b, e\}, d\}, c\}, a\} -$$ $$ \{\{\{\{c, e\}, b\}, d\}, a\} + 
    \{\{\{\{d, e\}, a\}, b\}, c\} - \{\{\{\{d, e\}, b\}, a\}, c\} + 
    \{\{\{\{d, e\}, b\}, c\}, a\},
$$

$$g_{15}(a,b,c,d,e)=$$
$$
 -\{\{\{a, b\}, c\}, \{d, e\}\} + 
    \{\{\{a, d\}, c\}, \{b, e\}\} + \{\{\{b, c\}, a\}, \{d, e\}\} + 
    \{\{\{b, c\}, e\}, \{a, d\}\} - $$ $$\{\{\{b, e\}, c\}, \{a, d\}\} - 
    \{\{\{c, d\}, a\}, \{b, e\}\} - \{\{\{c, d\}, e\}, \{a, b\}\} + 
    \{\{\{d, e\}, c\}, \{a, b\}\} +$$ $$ \{\{\{\{a, b\}, c\}, d\}, e\} - 
    \{\{\{\{a, d\}, c\}, b\}, e\} - \{\{\{\{b, c\}, a\}, d\}, e\} - 
    \{\{\{\{b, c\}, e\}, d\}, a\} +$$ $$ \{\{\{\{b, e\}, c\}, d\}, a\} + 
    \{\{\{\{c, d\}, a\}, b\}, e\} + \{\{\{\{c, d\}, e\}, b\}, a\} - 
    \{\{\{\{d, e\}, c\}, b\}, a\}, 
$$

$$g_{16}(a,b,c,d,e)=$$
$$
-\{\{\{a, b\}, c\}, \{d, e\}\} + 
    \{\{\{a, e\}, c\}, \{b, d\}\} + \{\{\{b, c\}, a\}, \{d, e\}\} + 
    \{\{\{b, c\}, d\}, \{a, e\}\} -$$ $$ \{\{\{b, d\}, c\}, \{a, e\}\} - 
    \{\{\{c, e\}, a\}, \{b, d\}\} - \{\{\{c, e\}, d\}, \{a, b\}\} + 
    \{\{\{d, e\}, c\}, \{a, b\}\} + $$ $$\{\{\{\{a, b\}, c\}, e\}, d\} - 
    \{\{\{\{a, e\}, c\}, b\}, d\} - \{\{\{\{b, c\}, a\}, e\}, d\} - 
    \{\{\{\{b, c\}, d\}, e\}, a\} + $$ $$\{\{\{\{b, d\}, c\}, e\}, a\} + 
    \{\{\{\{c, e\}, a\}, b\}, d\} + \{\{\{\{c, e\}, d\}, b\}, a\} - 
    \{\{\{\{d, e\}, c\}, b\}, a\},
$$

$$g_{17}(a,b,c,d,e)=$$ $$
 -\{\{\{a, b\}, e\}, \{c, d\}\} + 
    \{\{\{a, e\}, d\}, \{b, c\}\} - \{\{\{b, c\}, d\}, \{a, e\}\} + 
    \{\{\{b, d\}, c\}, \{a, e\}\} +$$ $$ \{\{\{b, e\}, a\}, \{c, d\}\} + 
    \{\{\{c, d\}, e\}, \{a, b\}\} - \{\{\{d, e\}, a\}, \{b, c\}\} - 
    \{\{\{d, e\}, c\}, \{a, b\}\} + $$ $$\{\{\{\{a, b\}, e\}, d\}, c\} - 
    \{\{\{\{a, e\}, d\}, b\}, c\} + \{\{\{\{b, c\}, d\}, e\}, a\} - 
    \{\{\{\{b, d\}, c\}, e\}, a\} -$$ $$ \{\{\{\{b, d\}, e\}, a\}, c\} + 
    \{\{\{\{b, d\}, e\}, c\}, a\} - \{\{\{\{b, e\}, a\}, d\}, c\} + 
    \{\{\{\{b, e\}, d\}, a\}, c\} -$$ $$ \{\{\{\{b, e\}, d\}, c\}, a\} - 
    \{\{\{\{c, d\}, e\}, b\}, a\} + \{\{\{\{d, e\}, a\}, b\}, c\} + 
    \{\{\{\{d, e\}, c\}, b\}, a\},
$$

$$g_{18}(a,b,c,d,e)=$$
$$
 -\{\{\{a, b\}, d\}, \{c, e\}\} + 
    \{\{\{a, e\}, d\}, \{b, c\}\} - \{\{\{b, c\}, d\}, \{a, e\}\} + 
    \{\{\{b, d\}, a\}, \{c, e\}\} +$$ $$ \{\{\{b, d\}, c\}, \{a, e\}\} + 
    \{\{\{c, e\}, d\}, \{a, b\}\} - \{\{\{d, e\}, a\}, \{b, c\}\} - 
    \{\{\{d, e\}, c\}, \{a, b\}\} +$$ $$ \{\{\{\{a, b\}, d\}, e\}, c\} - 
    \{\{\{\{a, e\}, d\}, b\}, c\} + \{\{\{\{b, c\}, d\}, e\}, a\} - 
    \{\{\{\{b, d\}, a\}, e\}, c\} - $$ $$\{\{\{\{b, d\}, c\}, e\}, a\} - 
    \{\{\{\{c, e\}, d\}, b\}, a\} + \{\{\{\{d, e\}, a\}, b\}, c\} + 
    \{\{\{\{d, e\}, c\}, b\}, a\},
$$

$$g_{19}(a,b,c,d,e)=$$
$$
 -\{\{\{a, e\}, c\}, \{b, d\}\} + 
    \{\{\{a, e\}, d\}, \{b, c\}\} - \{\{\{b, c\}, d\}, \{a, e\}\} + 
    \{\{\{b, d\}, c\}, \{a, e\}\} +$$ $$ \{\{\{c, e\}, a\}, \{b, d\}\} - 
    \{\{\{d, e\}, a\}, \{b, c\}\} + \{\{\{\{a, e\}, c\}, d\}, b\} - 
    \{\{\{\{a, e\}, d\}, c\}, b\} +$$ $$ \{\{\{\{b, c\}, d\}, e\}, a\} - 
    \{\{\{\{b, d\}, c\}, e\}, a\} - \{\{\{\{c, e\}, a\}, d\}, b\} + 
    \{\{\{\{c, e\}, d\}, a\}, b\} -$$ $$ \{\{\{\{c, e\}, d\}, b\}, a\} + 
    \{\{\{\{d, e\}, a\}, c\}, b\} - \{\{\{\{d, e\}, c\}, a\}, b\} + 
    \{\{\{\{d, e\}, c\}, b\}, a\},
$$

$$g_{20}(a,b,c,d,e)=$$
$$
 -\{\{\{a, d\}, c\}, \{b, e\}\} + 
    \{\{\{a, e\}, d\}, \{b, c\}\} - \{\{\{b, c\}, d\}, \{a, e\}\} + 
    \{\{\{b, e\}, c\}, \{a, d\}\} +$$ $$ \{\{\{c, d\}, a\}, \{b, e\}\} + 
    \{\{\{c, d\}, b\}, \{a, e\}\} - \{\{\{c, e\}, b\}, \{a, d\}\} - 
    \{\{\{d, e\}, a\}, \{b, c\}\} +$$ $$ \{\{\{\{a, d\}, c\}, e\}, b\} - 
    \{\{\{\{a, e\}, d\}, c\}, b\} + \{\{\{\{b, c\}, d\}, e\}, a\} - 
    \{\{\{\{b, e\}, c\}, d\}, a\} -$$ $$ \{\{\{\{c, d\}, a\}, e\}, b\} - 
    \{\{\{\{c, d\}, b\}, e\}, a\} + \{\{\{\{c, e\}, b\}, d\}, a\} + 
    \{\{\{\{c, e\}, d\}, a\}, b\} -$$ $$ \{\{\{\{c, e\}, d\}, b\}, a\} + 
    \{\{\{\{d, e\}, a\}, c\}, b\} - \{\{\{\{d, e\}, c\}, a\}, b\} + 
    \{\{\{\{d, e\}, c\}, b\}, a\}.
$$

Below we apply  the following re-writing rules for monomials of degree $5.$ 

$$a_j\circ a_i\rightarrow a_i\circ  a_j, \quad i<j,$$
$$(a_i\circ a_j)\circ (a_s\circ a_k)\rightarrow (a_s\circ a_k)\circ(a_i\circ a_j),\quad i>s, i<j,s<k, $$
$$(x\circ (y\circ z))\circ(u\circ v)\rightarrow ((y\circ z)\circ x)\circ(u\circ v),$$
$$((x\circ (y\circ z))\circ u)\circ v\rightarrow (((y\circ z)\circ x)\circ u)\circ v.$$
Here to simplify notations we set 
$a_1=a, a_2=b,a_3=c,a_4=d, a_5=e$
We endow generators by order, 
$$a_i<a_j \mbox{\;iff \;$i<j$ }$$
and we assume that $x,y,z,u,v\in \{a,b,c,d,e\}.$
Then we obtain the following presentations $g_j(a,b,c,d,e),1\le j\le 20,$
as a linear combinations of polynomials of a form  $f_5(a_{i_1},a_{i_2},a_{i_3},a_{i_4},a_{i_5}).$

\bigskip

$$3\,g_1(a,b,c,d,e)=$$ 
  $$f_5(a, c, b, d, e) - 
     f_5(a, d, b, c, e) - 
     2 f_5(a, d, c, b, e) + 
     f_5(a, e, b, c, d) + 
     2 f_5(a, e, c, b, d) - $$ $$
     2 f_5(a, e, d, b, c) + 
     2 f_5(b, c, a, d, e) - 
     2 f_5(b, d, a, c, e) + 
     2 f_5(b, d, c, a, e) + 
     2 f_5(b, e, a, c, d) - $$ $$
     2 f_5(b, e, c, a, d) + 
     2 f_5(b, e, d, a, c) - 
     f_5(c, d, a, b, e) + 
     f_5(c, d, b, a, e) + 
     f_5(c, e, a, b, d) - $$ $$
     f_5(c, e, b, a, d) - 
     2 f_5(c, e, d, a, b) - 
     f_5(d, e, a, b, c) + 
     f_5(d, e, b, a, c) - 
     f_5(d, e, c, a, b),$$

    $$3\,g_2(a,b,c,d,e)=$$ 
  $$-f_5(a, c, b, d, e) + 
     f_5(a, d, b, c, e) + 
     2 f_5(a, d, c, b, e) - 
     f_5(a, e, b, c, d) - $$ $$
     2 f_5(a, e, c, b, d) + 
     2 f_5(a, e, d, b, c) - 
     2 f_5(b, c, a, d, e) + 
     2 f_5(b, d, a, c, e) - $$ $$
     2 f_5(b, d, c, a, e) - 
     2 f_5(b, e, a, c, d) + 
     2 f_5(b, e, c, a, d) - 
     2 f_5(b, e, d, a, c) + $$ $$
     f_5(c, d, a, b, e) - 
     f_5(c, d, b, a, e) - 
     f_5(c, e, a, b, d) + 
     f_5(c, e, b, a, d) - $$ $$
     f_5(c, e, d, a, b) + 
     f_5(d, e, a, b, c) - 
     f_5(d, e, b, a, c) - 
     2 f_5(d, e, c, a, b),
     $$

 $$3\,g_3(a,b,c,d,e)=$$ 
  $$2 f_5(a, c, b, d, e) + 
     f_5(a, d, b, c, e) - 
     f_5(a, d, c, b, e) - 
     f_5(a, e, b, c, d) + 
     f_5(a, e, c, b, d) - $$ $$
     f_5(a, e, d, b, c) + 
     f_5(b, c, a, d, e) - 
     f_5(b, d, a, c, e) + 
     f_5(b, d, c, a, e) + 
     f_5(b, e, a, c, d) - $$ $$
     f_5(b, e, c, a, d) + 
     f_5(b, e, d, a, c) + 
     f_5(c, d, a, b, e) - 
     f_5(c, d, b, a, e) - 
     f_5(c, e, a, b, d) + $$ $$
     f_5(c, e, b, a, d) - 
     f_5(c, e, d, a, b) - 
     2 f_5(d, e, a, b, c) + 
     2 f_5(d, e, b, a, c) - 
     2 f_5(d, e, c, a, b),
     $$

$$3g_4(a,b,c,d,e)=$$ 

    $$-2 f_5(a, c, b, d, e) - 
     f_5(a, d, b, c, e) + 
     f_5(a, d, c, b, e) + 
     f_5(a, e, b, c, d) - 
     f_5(a, e, c, b, d) + $$
     $$
     f_5(a, e, d, b, c) - 
     f_5(b, c, a, d, e) + 
     f_5(b, d, a, c, e) - 
     f_5(b, d, c, a, e) - 
     f_5(b, e, a, c, d) + $$
     $$
     f_5(b, e, c, a, d) - 
     f_5(b, e, d, a, c) - 
     f_5(c, d, a, b, e) + 
     f_5(c, d, b, a, e) + 
     f_5(c, e, a, b, d) - $$ 
    $$
     f_5(c, e, b, a, d) + 
     f_5(c, e, d, a, b) - 
     f_5(d, e, a, b, c) + 
     f_5(d, e, b, a, c) - 
     f_5(d, e, c, a, b)
$$

$$3g_5(a,b,c,d,e)=$$ 
$$f_5(a, c, b, d, e) - 
     f_5(a, d, b, c, e) - 
     2 f_5(a, d, c, b, e) + 
     f_5(a, e, b, c, d) + 
     2 f_5(a, e, c, b, d) + $$ $$
     f_5(a, e, d, b, c) + 
     2 f_5(b, c, a, d, e) - 
     2 f_5(b, d, a, c, e) + 
     2 f_5(b, d, c, a, e) + 
     2 f_5(b, e, a, c, d) - $$ $$
     2 f_5(b, e, c, a, d) + 
     2 f_5(b, e, d, a, c) - 
     f_5(c, d, a, b, e) + 
     f_5(c, d, b, a, e) + 
     f_5(c, e, a, b, d) - $$ $$
     f_5(c, e, b, a, d) + 
     f_5(c, e, d, a, b) - 
     f_5(d, e, a, b, c) + 
     f_5(d, e, b, a, c) - 
     f_5(d, e, c, a, b))
$$

$$3g_6(a,b,c,d,e)=$$ 
  $$-f_5(a, c, b, d, e) - 
     2 f_5(a, d, b, c, e) - 
     f_5(a, d, c, b, e) + 
     2 f_5(a, e, b, c, d) + 
     f_5(a, e, c, b, d) + $$ $$
     2 f_5(a, e, d, b, c) + 
     f_5(b, c, a, d, e) - 
     f_5(b, d, a, c, e) + 
     f_5(b, d, c, a, e) + 
     f_5(b, e, a, c, d) - $$ $$
     f_5(b, e, c, a, d) - 
     2 f_5(b, e, d, a, c) - 
     2 f_5(c, d, a, b, e) + 
     2 f_5(c, d, b, a, e) + 
     2 f_5(c, e, a, b, d) - $$ $$
     2 f_5(c, e, b, a, d) + 
     2 f_5(c, e, d, a, b) + 
     f_5(d, e, a, b, c) - 
     f_5(d, e, b, a, c) + 
     f_5(d, e, c, a, b),$$

$$3g_7(a,b,c,d,e)=$$ $$
  -f_5(a, c, b, d, e) - 
     2 f_5(a, d, b, c, e) - 
     f_5(a, d, c, b, e) + 
     2 f_5(a, e, b, c, d) + 
     f_5(a, e, c, b, d) + $$ $$
     2 f_5(a, e, d, b, c) + 
     f_5(b, c, a, d, e) - 
     f_5(b, d, a, c, e) + 
     f_5(b, d, c, a, e) + 
     f_5(b, e, a, c, d) - $$ $$
     f_5(b, e, c, a, d) - 
     2 f_5(b, e, d, a, c) - 
     2 f_5(c, d, a, b, e) + 
     2 f_5(c, d, b, a, e) + 
     2 f_5(c, e, a, b, d) - $$ $$
     2 f_5(c, e, b, a, d) + 
     2 f_5(c, e, d, a, b) + 
     f_5(d, e, a, b, c) - 
     f_5(d, e, b, a, c) + 
     f_5(d, e, c, a, b),$$

 $$3g_8(a,b,c,d,e)=$$ 
  $$-f_5(a, c, b, d, e) + 
     f_5(a, d, b, c, e) + 
     2 f_5(a, d, c, b, e) - 
     f_5(a, e, b, c, d) - 
     2 f_5(a, e, c, b, d) - $$ $$
     f_5(a, e, d, b, c) - 
     2 f_5(b, c, a, d, e) + 
     2 f_5(b, d, a, c, e) - 
     2 f_5(b, d, c, a, e) - 
     2 f_5(b, e, a, c, d) + $$ $$
     2 f_5(b, e, c, a, d) - 
     2 f_5(b, e, d, a, c) + 
     f_5(c, d, a, b, e) - 
     f_5(c, d, b, a, e) - 
     f_5(c, e, a, b, d) + $$ $$
     f_5(c, e, b, a, d) - 
     f_5(c, e, d, a, b) - 
     2 f_5(d, e, a, b, c) - 
     f_5(d, e, b, a, c) - 
     2 f_5(d, e, c, a, b),$$

$$3g_{9}(a,b,c,d,e)=$$ 
$$
2 f_5(a, c, b, d, e) + 
     f_5(a, d, b, c, e) - 
     f_5(a, d, c, b, e) + 
     2 f_5(a, e, b, c, d) + 
     f_5(a, e, c, b, d) - $$ $$
     f_5(a, e, d, b, c) + 
     f_5(b, c, a, d, e) - 
     f_5(b, d, a, c, e) + 
     f_5(b, d, c, a, e) + 
     f_5(b, e, a, c, d) - $$ $$
     f_5(b, e, c, a, d) + 
     f_5(b, e, d, a, c) + 
     f_5(c, d, a, b, e) - 
     f_5(c, d, b, a, e) - 
     f_5(c, e, a, b, d) + $$ $$
     f_5(c, e, b, a, d) - 
     f_5(c, e, d, a, b) - 
     2 f_5(d, e, a, b, c) + 
     2 f_5(d, e, b, a, c) - 
     2 f_5(d, e, c, a, b),$$

$$3g_{10}(a,b,c,d,e)=$$ 
$$-2 f_5(a, c, b, d, e) + 
     2 f_5(a, d, b, c, e) + 
     f_5(a, d, c, b, e) + 
     f_5(a, e, b, c, d) - 
     f_5(a, e, c, b, d) + $$ $$
     f_5(a, e, d, b, c) - 
     f_5(b, c, a, d, e) + 
     f_5(b, d, a, c, e) - 
     f_5(b, d, c, a, e) - 
     f_5(b, e, a, c, d) + $$ $$
     f_5(b, e, c, a, d) - 
     f_5(b, e, d, a, c) - 
     f_5(c, d, a, b, e) + 
     f_5(c, d, b, a, e) + 
     f_5(c, e, a, b, d) - $$ $$
     f_5(c, e, b, a, d) + 
     f_5(c, e, d, a, b) - 
     f_5(d, e, a, b, c) + 
     f_5(d, e, b, a, c) - 
     f_5(d, e, c, a, b),$$

 $$g_{11}(a,b,c,d,e)=-f_5(a, d, b, c, e) - f_5(b, d, a, c, e),$$

 $$3g_{12}(a,b,c,d,e)=$$ 
 $$
-f_5(a, c, b, d, e) - 
     2 f_5(a, d, b, c, e) - 
     f_5(a, d, c, b, e) - 
     f_5(a, e, b, c, d) + 
     f_5(a, e, c, b, d) + $$ $$
     2 f_5(a, e, d, b, c) + 
     f_5(b, c, a, d, e) - 
     f_5(b, d, a, c, e) + 
     f_5(b, d, c, a, e) - 
     2 f_5(b, e, a, c, d) - $$ $$
     f_5(b, e, c, a, d) - 
     2 f_5(b, e, d, a, c) - 
     2 f_5(c, d, a, b, e) + 
     2 f_5(c, d, b, a, e) + 
     2 f_5(c, e, a, b, d) - $$ $$
     2 f_5(c, e, b, a, d) + 
     2 f_5(c, e, d, a, b) + 
     f_5(d, e, a, b, c) - 
     f_5(d, e, b, a, c) + 
     f_5(d, e, c, a, b),
     $$

     $$3g_{13}(a,b,c,d,e)=$$ 
  $$-f_5(a, c, b, d, e) - 
     2 f_5(a, d, b, c, e) - 
     f_5(a, d, c, b, e) - 
     f_5(a, e, b, c, d) + 
     f_5(a, e, c, b, d) - $$ $$
     f_5(a, e, d, b, c) + 
     f_5(b, c, a, d, e) - 
     f_5(b, d, a, c, e) + 
     f_5(b, d, c, a, e) + 
     f_5(b, e, a, c, d) - $$ $$
     f_5(b, e, c, a, d) + 
     f_5(b, e, d, a, c) - 
     2 f_5(c, d, a, b, e) + 
     2 f_5(c, d, b, a, e) - 
     f_5(c, e, a, b, d) + $$ $$
     f_5(c, e, b, a, d) - 
     f_5(c, e, d, a, b) + 
     f_5(d, e, a, b, c) - 
     f_5(d, e, b, a, c) + 
     f_5(d, e, c, a, b), $$

$$3g_{14}(a,b,c,d,e)=$$ 
    $$2 f_5(a, c, b, d, e) - 
     2 f_5(a, d, b, c, e) - 
     f_5(a, d, c, b, e) - 
     f_5(a, e, b, c, d) + 
     f_5(a, e, c, b, d) - $$ $$
     f_5(a, e, d, b, c) + 
     f_5(b, c, a, d, e) - 
     f_5(b, d, a, c, e) + 
     f_5(b, d, c, a, e) + 
     f_5(b, e, a, c, d) - $$ $$
     f_5(b, e, c, a, d) + 
     f_5(b, e, d, a, c) - 
     2 f_5(c, d, a, b, e) + 
     2 f_5(c, d, b, a, e) - 
     f_5(c, e, a, b, d) + $$ $$
     f_5(c, e, b, a, d) - 
     f_5(c, e, d, a, b) + 
     f_5(d, e, a, b, c) - 
     f_5(d, e, b, a, c) + 
     f_5(d, e, c, a, b),$$

     $$g_{15}(a,b,c,d,e)= f_5(a, e, c, b, d),$$

$$g_{16}(a,b,c,d,e)=f_5(a, d, c, b, e),$$

$$3g_{17}(a,b,c,d,e)=$$ 
$$f_5(a, c, b, d, e) - 
     f_5(a, d, b, c, e) - 
     2 f_5(a, d, c, b, e) + 
     f_5(a, e, b, c, d) - 
     f_5(a, e, c, b, d) - $$ $$
     2 f_5(a, e, d, b, c) + 
     2 f_5(b, c, a, d, e) - 
     2 f_5(b, d, a, c, e) + 
     2 f_5(b, d, c, a, e) + 
     2 f_5(b, e, a, c, d) - $$ $$
     2 f_5(b, e, a, c, d) - 
     2 f_5(b, e, c, a, d) + 
     2 f_5(b, e, d, a, c) - 
     f_5(c, d, a, b, e) + 
     f_5(c, d, b, a, e) - $$ $$
     2 f_5(c, e, a, b, d) - 
     f_5(c, e, b, a, d) - 
     2 f_5(c, e, d, a, b) - 
     f_5(d, e, a, b, c) + 
     f_5(d, e, b, a, c) - 
     f_5(d, e, c, a, b),$$

$$g_{18}(a,b,c,d,e)= -f_5(a, d, c, b, e) - f_5(c, d, a, b, e),$$

   $$3g_{19}(a,b,c,d,e)=$$ 
  $$f_5(a, c, b, d, e) - 
     f_5(a, d, b, c, e) - 
     2 f_5(a, d, c, b, e) + 
     f_5(a, e, b, c, d) - 
     f_5(a, e, c, b, d) + $$ $$
     f_5(a, e, d, b, c) + 
     2 f_5(b, c, a, d, e) - 
     2 f_5(b, d, a, c, e) + 
     2 f_5(b, d, c, a, e) - 
     f_5(b, e, a, c, d) + $$ $$
     f_5(b, e, c, a, d) - 
     f_5(b, e, d, a, c) - 
     f_5(c, d, a, b, e) + 
     f_5(c, d, b, a, e) + 
     f_5(c, e, a, b, d) - $$ $$
     f_5(c, e, b, a, d) + 
     f_5(c, e, d, a, b) - 
     f_5(d, e, a, b, c) + 
     f_5(d, e, b, a, c) - 
     f_5(d, e, c, a, b)),$$

$$3g_{20}(a,b,c,d,e)=$$ 
$$  -2 f_5(a, c, b, d, e) - 
     f_5(a, d, b, c, e) - 
     2 f_5(a, d, c, b, e) + 
     f_5(a, e, b, c, d) - 
     f_5(a, e, c, b, d) + $$ $$
     f_5(a, e, d, b, c) - 
     f_5(b, c, a, d, e) - 
     2 f_5(b, d, a, c, e) + 
     2 f_5(b, d, c, a, e) - 
     f_5(b, e, a, c, d) + $$ $$
     f_5(b, e, c, a, d) - 
     f_5(b, e, d, a, c) - 
     f_5(c, d, a, b, e) + 
     f_5(c, d, b, a, e) + 
     f_5(c, e, a, b, d) - $$ $$
     f_5(c, e, b, a, d) + 
     f_5(c, e, d, a, b) - 
     f_5(d, e, a, b, c) + 
     f_5(d, e, b, a, c) - 
     f_5(d, e, c, a, b)).$$ 

     So, we have established that any multi-linear identity of degree $5$ for the algebra $AR^+$ is a consequence of the identity $f_5(a,b,c,d,e)=0.$ $\square$

\section{Algebras with multiplication $m_\eps$}

Let
$$m_\eps(a,b)=a\int\!\!\int b+\eps\int a \int b$$
and
$$(a,b,c)_\eps=m_\eps (a,m_\eps(b,c))-m_\eps(m_\eps(a,b),c)$$
We have 
$$jor(m_1)(x^i,x^j)=\frac{(i+j+3) (i+j+4) x^{i+j+2}}{(i+1) (i+2) (j+1) (j+2)}$$
Therefore,
$$jor(m_1)(jor(m_1)(x^i,x^j),x^s)=
\frac{(5 + i + j + s) (6 + i + j + s) x^{
  4 + i + j + s}}{
(1 + i) (2 + i) (1 + j) (2 + j) (1 + s) (2 + s)}
$$
and
$$jor(m_1)(x^i,jor(m_1)(x^j,x^s))=
\frac{(5 + i + j + s) (6 + i + j + s) x^{
  4 + i + j + s}}{
(1 + i) (2 + i) (1 + j) (2 + j) (1 + s) (2 + s)}
$$
So,
$$(a,b,c)_1=0$$
is identity.

\sudda{
Symmetries of $baxter2$ in degree $5:$

$$bbaxter2[n][n[a,e],b,c,d]-n[bbaxter2[n][a,b,c,d],e]=0,$$

$$
bbaxter2[n][{a,b,n[c,d],e}]-bbaxter2[n][{a,b,n[c,e],d}]+$$ $$bbaxter2[n][{a,c,n[b,d],e}]-bbaxter2[n][{a,c,n[b,e],d}]=0,
$$

$$
n[a,bbaxter2[n][{b,c,d,e}]]- n[a,bbaxter2[n][{c,b,d,e}]]+
n[a,bbaxter2[n][{d,b,c,e}]]$$
$$
-n[bbaxter2[n][{a,b,c,e}],d] + n[bbaxter2[n][{a,b,d,e}],c] - 
 n[bbaxter2[n][{a,c,d,e}],b]$$
$$
-bbaxter2[n][{a,b,c,n[d,e]}] + bbaxter2[n][{a,b,d,n[c,e]}] - 
 bbaxter2[n][{a,c,d,n[b,e]}]$$
$$
+bbaxter2[n][{a,b,n[c,e],d}] - bbaxter2[n][{a,b,n[d,e],c}] - 
 bbaxter2[n][{a,c,n[b,d],e}] + $$ $$bbaxter2[n][{a,c,n[d,e],b}] + 
 bbaxter2[n][{a,d,n[b,c],e}] - bbaxter2[n][{a,d,n[c,b],e}]=0.
 $$
 
All of these identities follows from the identity $rcom=0.$
Therefore the identity $f_5^+(a,b,c,d,e)=0$ is not consequence of these identities, since the identity $f_5^+=0$ is not consequence of identity $rcom=0.$
}

\section{An algebra with multiplication $e_i\circ e_j=\eps_{i,j} e_i,$ 
where $\eps_{i,j}=0,$ if $i\le j$ and $\eps_{i,j}=\eps_j,$ if $i>j.$}


Let $M=(\eps_{i,j})_{i,j\in I}$
be matrix whose components are indexed by elements of some subset $I\subseteq {\bf Z}$
such that 
$$\eps_{i,j}=0, \quad i\le j.$$
$$\eps_{i,j}=\eps_j, \quad i>j,$$
where $\eps_j\in K.$
For example, if $I=\{1,2,\ldots,n\},$ then 
$$M=\left(\begin{array}{ccccccc}
0&0&0&\cdots &0&0\\
\eps_{1}&0&0&\cdots &0&0\\
\eps_{1}&\eps_2&0&\cdots &0&0\\
\eps_{1}&\eps_2&\eps_3&\cdots &0&0\\
\vdots&\vdots&\vdots&\cdots &\vdots&\vdots\\
\eps_1&\eps_2&\eps_3&\cdots&0&0\\
\eps_1&\eps_2&\eps_3&\cdots&\eps_{n-1}&0
\end{array}\right)
$$

Since
$$(e_i\circ e_j)\circ e_s=
(\eps_{i,j} e_i)\circ e_s=\eps_{i,j}\eps_{i,s} e_i,
$$
$$(e_i\circ e_s)\circ e_j=
(\eps_{i,s} e_i)\circ e_j=\eps_{i,s}\eps_{i,j} e_i,
$$
right-commutative identity 
$$(e_i\circ e_j)\circ e_s=(e_i\circ e_s)\circ e_j.$$
holds for any $i,j\in I.$

Let us prove that
$$f_4(e_i,e_j,e_s,e_k)=0$$
is identity for any $i,j,s,k\in I.$
We have 
$$e_i\circ ([e_j,e_s]\circ e_k)=$$

$$\eps_{j,s}e_i\circ (e_j\circ e_k)-
\eps_{s,j}e_i\circ (e_s\circ e_k)=$$

$$\eps_{j,s}\eps_{j,k}e_i\circ e_j-
\eps_{s,j}\eps_{s,k}e_i\circ e_s=$$

$$\eps_{j,s}\eps_{j,k}\eps_{i,j}e_i-
\eps_{s,j}\eps_{s,k}\eps_{i,s} e_i=$$

$$(\eps_{i,j}\eps_{j,s}\eps_{j,k}-
\eps_{i,s}\eps_{s,j}\eps_{s,k}) e_i,$$
Further,
$$-(e_i,e_j,e_s\circ e_k)=$$

$$-e_i\circ (e_j\circ (e_s\circ e_k))+(e_i\circ e_j)\circ (e_s\circ e_k)=$$

$$-\eps_{s,k} e_i\circ (e_j\circ e_s)+\eps_{i,j}\eps_{s,k}\, e_i\circ e_s=$$

$$-\eps_{s,k}\eps_{j,s} \eps_{i,j}e_i+\eps_{i,j}\eps_{s,k}\eps_{i,s}\, e_i=$$

$$(- \eps_{i,j}\eps_{j,s}\eps_{s,k}
+\eps_{i,j}\eps_{i,s}\eps_{s,k})e_i,$$
Similarly,
$$+(e_i,e_s,e_j\circ e_k)=$$

$$(- \eps_{i,s}\eps_{s,j}\eps_{j,k}
+\eps_{i,s}\eps_{i,j}\eps_{j,k})e_i.$$
Hence,
$$f_4(e_i,e_j,e_s,e_k)=$$

$$(\eps_{i,j}\eps_{j,s}\eps_{j,k}-
\eps_{i,s}\eps_{s,j}\eps_{s,k}) e_i+$$
$$(- \eps_{i,j}\eps_{j,s}\eps_{s,k}
+\eps_{i,j}\eps_{i,s}\eps_{s,k})e_i+$$
$$(- \eps_{i,s}\eps_{s,j}\eps_{j,k}
+\eps_{i,s}\eps_{i,j}\eps_{j,k})e_i=$$

\sudda{
$$\{\eps_{i,j}\eps_{j,s}\eps_{j,k}-
\eps_{i,s}\eps_{s,j}\eps_{s,k}- \eps_{i,j}\eps_{j,s}\eps_{s,k}
+\eps_{i,j}\eps_{i,s}\eps_{s,k}- \eps_{i,s}\eps_{s,j}\eps_{j,k}
+\eps_{i,s}\eps_{i,j}\eps_{j,k}\}e_i=$$
}
$$G_{i,j,s,k}e_i,$$
where
$$G_{i,j,s,k}=$$
$$\eps_{i,j}\eps_{j,s}\eps_{j,k}
+\eps_{i,j}\eps_{i,s}\eps_{s,k}
+\eps_{i,s}\eps_{i,j}\eps_{j,k}
-\eps_{i,s}\eps_{s,j}\eps_{s,k}
- \eps_{i,j}\eps_{j,s}\eps_{s,k}
- \eps_{i,s}\eps_{s,j}\eps_{j,k}
.$$

We have to establish that
$$G_{i,j,s,k}=0$$
for any $i,j,s,k\in I. $
Since,
$$G_{i,j,s,k}+G_{i,s,j,k}=0,$$
it is enough to do that for the case $j<s.$

Might be the following three cases.

\bigskip

Case I. Let $i<j.$ Then $\eps_{i,j}=0,$ $\eps_{i,s}=0.$ 

\bigskip

Case II. Let $j\le i \le s.$ Then $\eps_{j,s}=\eps_{i,s}=0.$ Hence
$$\eps_{i,j}\; \eps_{j,s}\;\eps_{j,k}=\eps_{i,j}\; 0\; \eps_{j,k}=0,$$
$$+\eps_{i,j}\;\eps_{i,s}\;\eps_{s,k}=
+\eps_{i,j}\;0\;\eps_{s,k}=0,$$
$$+\eps_{i,s}\;\eps_{i,j}\eps_{j,k}=0\;\eps_{i,j}\eps_{j,k}=0,$$
$$-\eps_{i,s}\;\eps_{s,j}\eps_{s,k}=
-0\;\eps_{s,j}\eps_{s,k}=0,$$
$$- \eps_{i,j}\;\eps_{j,s}\;\eps_{s,k}
- \eps_{i,j}\;0\;\eps_{s,k}=0,$$
$$- \eps_{i,s}\;\eps_{s,j}\eps_{j,k}=
-0\;\eps_{s,j}\eps_{j,k}=0.$$

\bigskip

Case III. Let $j<s<i.$ Then $\eps_{j,s}=0, \eps_{s,i}=0,$
and $\eps_{i,j}=\eps_j, \eps_{i,s}=\eps_s, \eps_{s,j}=\eps_j.$
Thus,
$$\eps_{i,j}\;\eps_{j,s}\;\eps_{j,k}=
\eps_{i,j}\;0\;\eps_{j,k}=0,$$

$$+\eps_{i,j}\eps_{i,s}\eps_{s,k}
-\eps_{i,s}\eps_{s,j}\eps_{s,k}=
\eps_j\eps_s\eps_{s,k}-\eps_s\eps j\eps_{s,k}=0,
$$

$$
+\eps_{i,s}\eps_{i,j}\eps_{j,k}
- \eps_{i,s}\eps_{s,j}\eps_{j,k}=
+\eps_{s}\eps_{j}\eps_{j,k}- \eps_{s}\eps_{j}\eps_{j,k}=0,
$$

$$- \eps_{i,j}\;\eps_{j,s}\;\eps_{s,k}=
- \eps_{i,j}\; 0\;\eps_{s,k}=0.$$
So, in all cases 
$$G_{i,j,s,k}=0.$$
So, we have established that 
$f_4(e_i,e_j,e_s,e_k)=0$ is identity. 

We see that
$$[e_i,e_j]=\pm \eps_{Min\{i,j\}}e_{Max\{i,j\}},$$
$$\{e_i,e_j\}=\eps_{Min\{i,j\}}e_{Max\{i,j\}}.$$
More exactly,
$$[e_i,e_j]=\left\{\begin{array}{cc}
\eps_j e_i,&\mbox{ if $i>j$,}\\
-\eps_i e_j,&\mbox{ if $i<j$,}
\end{array}\right.$$

$$\{e_i,e_j\}=\left\{\begin{array}{cc}
\eps_j e_i,&\mbox{ if $i>j$,}\\
2\eps_i e_i,&\mbox{ if $i=j$},\\
\eps_i e_j,&\mbox{ if $i<j$}.\\
\end{array}\right.$$

It seems that any such algebra under Lie commutator has ideal $J=\{ e_i, | i\in I, i\ne Min(I)\}$
since $e_{Min}(I)$ never appear in $A.$ 

If $I={\bf Z},$ then our algebra coincides with itself. Hence it is not solvable or nilpotent. 

\sudda{Is it our algebra associative-special?
Suppose that there exists associative multiplication $a\cdot b$ on a space $A=\{e_i | i\in I\}$ and linear map $R:A\rightarrow A,$ such that 
$$e_i\circ e_j=e_i \cdot R(e_j)$$
Suppose that
$$R(e_1)=\sum_{i\ge 1} \mu_i e_i$$
Then 
$$R(e_1)\cdot R(e_1)=R(e_1)\circ e_1=\eps_1\sum_{i\ge 2}\mu_i e_i,$$
and
$$(R(e_1)\cdot R(e_1))\cdot R(e_1)=(R(e_1)\circ e_1)\circ e_1=
\eps_1^2\sum_{i\ge 2}\mu_i^2 e_i,$$

$$R(e_1)\cdot (R(e_1)\cdot R(e_1))=
R(e_1)\circ (e_1\cdot R(e_1)+R(e_1)\cdot e_1+\lam e_1\cdot e_1)=??
$$
}


\section{The multiplication $\star_{k,n}$}
Let us define multiplication $a\star_{k,n} b$ by 
$$a\star_{k,n} b=\sum_{i=0}^n {n\choose i} \int_i a \int_{n-i+k}b$$
For example,
$$x^i\star_{0,0}x^j=x^{i+j},$$
$$x^i\star_{0,1}x^j=
\frac{(i+j+2) x^{i+j+1}}{(i+1)
   (j+1)},$$
$$x^i\star_{0,2}x^j=
\frac{(i+j+3) (i+j+4)
   x^{i+j+2}}{(i+1) (i+2) (j+1)
   (j+2)},$$

$$x^i\star_{1,0}x^j=\frac{x^{i+j+1}}{j+1},$$
$$x^i\star_{1,1}x^j=\frac{(i+j+3) x^{i+j+2}}{(i+1)
   (j+1) (j+2)},$$
   $$x^i\star_{1,2} x^j=\frac{(i+j+4) (i+j+5)
   x^{i+j+3}}{(i+1) (i+2) (j+1)
   (j+2) (j+3)}.$$

The multiplication $a\star_{k,n} b$ has the following properties 
\begin{itemize}
    \item $\star_{0,n}$ is associative and commutative (see section \ref{associative})
    \item $\star_{1,n}$ is not associative, not commutative, but left-zinbiel
    \item $\star_{2,n}$ is not left-zinbiel, does not satisfy the identity $f_4=0,$ but satisfy the identities $rcom=0$ and $f_4'=0.$
    \item $\star_{k,n}$ is right-commutative, if $k>2.$
\end{itemize}

Therefore, multiplications $\star_{0,n}$ 
satisfy the identities $rcom=0, f_4=0, f_4'=0, tortkara^{-}=0$ and $f_5^{+}=0.$    
The multiplication $\star_{2,n}$ 
satisfy identities $rcom=0, f_4'=0$ and 
$tortkara^-=0,$ but does not satisfy identities $f_4=0, f_5^+=0.$ 


Let us check  right-commutative property of $\star_{k,n}.$
We have 
$$(x^i\star_{k,n} x^j)\star_{k,n}x^s=\frac{{i+j+s+2k+3n\choose n}  x^{2(k+n)+i+j+s}}
{((k+n)!^2\;{i+k\choose n}
{s+k+n\choose k+n}{j+k+n\choose k+n}} $$
and for any $i,j,s\ge 0$
$$(x^i\star_{k,n}x^j)\star_{k,n}x^s=(x^i\star_{k,n}x^s)\star_{k,n}x^j.$$
So, the multiplication $\star_{k,n}$ is right-symmetric for any $k\ge 0.$

\section{Lie commutators of Zinbiel products}

Let 
$$[a,b]_{n}=\sum_{i=0}^n {n\choose i}\left(1-\frac{2i}{n}\right)\int_i a\int_{n-i}b,$$
where
$$\int_i a=\int\cdots \int a$$
Here number of integrals is $i$ and $\int a(x)=\int_0^x a(x)dx.$
 For example,
 $$[a,b]_1=a\int b-b\int a,$$
 
$$[a,b]_2=a\int_2 b-\int_2 a\; b,$$

$$[a,b]_3=a\int_3 b
+\int a\;\int_2 b
-\int_2 a\int b
-\int_3a\;b,$$

$$[a,b]_4=
a\int_4 b
+2\int a\int_3 b
-2\int_3 a\int b
-\int_4a \; b,$$

\sudda{
$$[a,b]_5=
a\int_5 b
+3\int a\int_4 b+
2\int_2 a\int_3 b
-2\int_3 a\int_2 b
-3\int_4 a\int b
-\int_5 a \; b,$$

$$[a,b]_6=
a\int_6 b
+4\int a\int_5 b+
5\int_2 a\int_4 b
-5\int_4 a\int_2 b
-4\int_5 a\int b
-\int_6 a \; b,$$}

and
$$[x^3,x^4]_1=-\frac{x^8}{20},\quad [x^3,x^4]_2=
-\frac{x^9}{60},\quad 
[x^3,x^4]_3=-\frac{11}{2100}x^{10}.$$

Recall that for any left-Zinbiel algebra $A$ with multiplication $ab$ its Lie commutator 
$[a,b]=ab-ba$ satisfies so called Tortkara identity,
$$[[a,b],[c,d]]+[[a,d],[c,b]]=[jac[a,b,c],d]+[jac[a,d,c],b],$$
where
$$jac[a,b,c]=[[a,b],c]+[[b,c],a]+[[c,a],b]$$
is Jacobian.

\begin{thm} 
For any $n\ge 1$ the algebra $(K[x], [\; ,\;]_n)$ is Tortkara.
\end{thm}

{\bf Proof.} It follows from the fact that 
$$[a,b]_{n+1}=a\star_{1,n} b-b\star_{1,n}a.$$
is a Lie commutator of left-Zinbiel multiplication $\star_{1,n}.$

\section{Multiplications with identity $f_4'(a,b,c,d)=0$ }

Let 
$$a\circ_2 b=\int a \int b+a \int_2 b,$$
$$a\circ_3 b=\int a \int_2b+a \int_3 b,$$
$$a\circ_4 b=\int_2 a \int_2 b+2 \int a \int_3 b+a \int_4 b.$$
Then
$$f_4'(a,b,c,d)=0$$
is identity for multiplications $a\circ_2 b,$ $a\circ_3b$ and $a\circ_4 b,$
but 
$$f_4(a,b,c,d)=0$$
is identity only for the multiplication $a\circ_2 b.$ Moreover,
the multiplication $a\circ_2 b$ satisfies left  Zinbiel identity 
$$(a\circ_2 b)\circ_2 c= a\circ_2(b\circ_2 c+c\circ_2 b).$$


\section{\label{associative}
The multiplication $a\star_{0,n} b=\sum_{i=0}^n {n\choose i}\int_i a \int_{n-i} b$ is associative}

Induction by $n.$
For $n=1$ this fact is well  known. 

Suppose that the multiplication $s\star_{0,n}$ is associative. Let us check that the map
$$\der: (A,\star_{0,n})\rightarrow (A,\star_{0,n+1})$$
is homomorphism.
We have 
$$\der(a\star_{0,n} b)=
\sum_{i=0}^n \der(\int_i a \int_{n-i} b)=$$

$$\der(a) \int_{n}b+a\int_{n-1}b+
\sum_{i=1}^n {n\choose i}\int_{i-1} a \int_{n-i} b+
\sum_{i=1}^{n-1} {n\choose i} \int_{i} a \int_{n-i-1} b
=$$

$$\der(a) \int_{n}b+ 
\sum_{i=0}^{n-1} {n+1\choose i+1}\int_{i} a \int_{n-i-1} b
=$$

$$\der(a) \int_{n+1}\der(b)+ 
\sum_{i=0}^{n-1} {n+1\choose i+1}\int_{i+1} \der(a) \int_{n-i} \der(b)
=$$

$$\der(a) \int_{n+1}\der(b)+ 
\sum_{i=1}^{n} {n+1\choose i}\int_{i} \der(a) \int_{n-i+1} \der(b)
=$$

$$\sum_{i=0}^{n+1} {n+1\choose i}\int_{i} \der(a) \int_{n-i+1} \der(b)=
$$

$$\der(a)\star_{0,n+1}\der(b)$$
Since any element of $A$ can be presented 
in a form $\der(a),$ 
from associativity of 
$(A,\star_{0,n})$ it follows that the algebra 
$(A,\star_{0,n+1})$ is also associative.

\section{$AR$-algebra of Novikov algebra}

Let $A$ be polynomial algebra $K[x].$ Take $R=\int$ and 
endow right-Novikov algebra $A$ by $AR$-multiplication 
$$a\diamond b=\der(a)\int b$$.

\begin{prp} \label{rsym} The algebra $(A,\diamond)$ is right-symmetric,
$$(a\diamond b)\diamond c-(a\diamond c)\diamond b=a\diamond (b\diamond c-c\diamond b),$$
where
$$[a,b]=a\diamond b-b\diamond a. $$
\end{prp}

{\bf Proof.}
By integration by parts we  have 
$$a\diamond(b\diamond c)=\der(a)\int( \der(b)\int c)=$$
$$\der (a)(b \int c-\int (b c)).$$
Further,
$$(a\diamond b)\diamond c=$$
$$\der(\der(a)\int b)\int c=$$

$$\der^2(a)\int b\int c+\der(a)b \int c.$$
Therefore,
$$(a,b,c)\stackrel{def}{=}a\diamond (b\diamond c)-(a\diamond b)\diamond c=$$

$$\der (a) (b \int c-\int (b c))-
\der^2(a)\int b\int c-\der(a)b \int c=$$

$$ -\der a\int(b c) -\der^2(a)\int b\int c.$$
Similarly, 
$$(a,c,b)=$$
$$ -\der a\int(c b) -\der^2(a)\int c\int b.$$
Since the algebra $A=K[x]$ is commutative, 
$$rsym(a,b,c)=(a,b,c)-(a,c,b)=0$$
is identity for the multiplication $\diamond$
$\square$

\begin{prp} Let 
$$s_{1,3}(t_1,t_2,t_3,t_4)=\sum_{\sigma\in S_4, \sigma(1)=1} sign\;\sigma\; ((t_{\sigma(1)}t_{\sigma(2)})t_{\sigma(3)})t_{\sigma(4)}.
$$
Then
$$s_{1,3}(a,b,c,d)=0$$
is identity of the algebra $(A,\diamond).$
\end{prp}

{\bf Proof.} Let 
$$S(a,b,c,d)=(((a\diamond b)\diamond c)\diamond d+
(((a\diamond c)\diamond d)\diamond b+
(((a\diamond d)\diamond b)\diamond c.$$

We have
$$(((a\diamond b)\diamond c)\diamond d=$$
$$ \der(\der(\der(a)\int b)\int c)\int d=$$
$$ \der(\der^2(a)\int b\int c+\der(a)\,b\int c)\int d=$$
$$ \der^3(a)\int b\int c\int d+
2\,\der^2(a) b\int c\int d
+\der^2(a)\int b\;c\int d+$$ $$
\der(a)\,\der(b)\int c\int d+
\der(a)\,b\, c\int d.$$
Similarly,
$$(((a\diamond c)\diamond d)\diamond b=$$
$$ \der^3(a)\int c\int d\int b+
2\der^2(a) c\int d\int b
+\der^2(a)\int c\;d\int b+$$ $$
\der(a)\,\der(c)\int d \int b+
\der(a)\,c\, d\int b,$$

$$(((a\diamond d)\diamond b)\diamond c=$$
$$ \der^3(a)\int d\int b\int c+
2\der^2(a) d\int b\int c
+\der^2(a)\int d\;b\int c+$$ $$
\der(a)\,\der(d)\int b \int c+
\der(a)\,d\, b\int c.$$
Thus,
$$S(a,b,c,d)=$$
$$(((a\diamond b)\diamond c)\diamond d+
(((a\diamond c)\diamond d)\diamond b+
(((a\diamond d)\diamond b)\diamond c=$$

$$ 3\der^3(a)\int b\int c\int d+$$
$$
3\,\der^2(a) (b\int c\int d+ c\int b\int d+ d\int b\int c)+$$
$$
\der(a)(\der(b)\int c\int d+
\der(c)\int d \int b+
\der(d)\int b \int c)+$$
$$\der(a)(b\, c\int d+b\,d\,\int c+ c\, d\,\int b).$$
So, $S(a,b,c,d)$ is symmetric under variables $b,c,d,$ and
$$s_{1,3}(a,b,c,d)=S(a,b,c,d)-S(a,c,b,d)=0.$$
$\square$

\begin{prp} \label{xxxx}The multiplication $\diamond$
satisfies the following  identity of degree $5$

$$(a\diamond  (b\diamond  c))\diamond [d, e]
- b\diamond  ((a\diamond [d,  e])\diamond  c)   
+ b\diamond  (d\diamond [c, a\diamond  e])
+ b\diamond  ((a\diamond  e)\diamond  [c, d]) $$ 
$$
+b\diamond  ([d,a\diamond  e]\diamond  c)  
   - b\diamond  ((a\diamond  d)\diamond [c, e])
   - b\diamond  (e\diamond  [c,  a\diamond  d])
   - b\diamond ( [e,  a\diamond  d]\diamond  c)$$      $$
   +d\diamond  [a\diamond  e,  b\diamond  c]
   - e\diamond  [a\diamond  d,b\diamond  c]
   + (a\diamond  [e,b\diamond  c])\diamond  d
 - (a\diamond  [d,  b\diamond  c])\diamond  e=0.$$
\end{prp}
    
{\bf Proof.}
It is enough to  check this fact for elements $a=e_{i_1},b=e_{i_2},c=e_{i_3},d=e_{i_4},e=e_{i_5},$ where 
$$e_i=x^{(i)}\stackrel{def}{=}\frac{x^i}{i!}$$
Then
$$e_i\diamond e_j={i+j\choose i-1}\;e_{i+j}.$$
$$[e_i,e_j]_\diamond=e_i\diamond e_j-e_j\diamond e_i={i+j+2\choose i+1}\frac{i-j}{i+j+2} \; e_{i+j},$$
$$
\{e_i,e_j\}_\diamond=e_i\diamond e_j+e_j\diamond e_i={i+j+2\choose i+1}\;\frac{i+j+i^2+j^2}{(i+j+1)(i+j+2)} \; e_{i+j}.$$
The calculation is straightforward. 

{\bf Remark.} Any Novikov algebra under Lie commutator satisfies standard skew-symmetric identity of degree $5$
$$\sum_{\sigma\in S_5,\sigma(1)=1} sign\,\sigma\, (((t_{\sigma(1)}t_{\sigma(2)} )t_{\sigma(3)})t_{\sigma(4)})t_{\sigma(5)}=0.$$
One can check that Novikov agebra satisfies the identity of degree $5$ given in Proposition  \ref{xxxx}.
But the multiplication $\diamond$ does not satisfy the standard skew-symmetric identity of degree $5.$


Let 
$$a\star b=ab-a\diamond b=a b-\der(a)\int b.$$
where $a,b\in K[x],$ such that $b(0)=0.$
It saitisfies the following identities
\begin{itemize}
    \item $rsym(a,b,c)=0$ is identity for the multiplication $a\star b=ab-\der(a)\int(b)$
    \item 
$$-2 \{a,c\}\star  (b\star  d) + 2 \{a,  d\}\star  (b\star  c) + 
  2 \{b, c\}\star  (a\star  d) - 2 \{b,d\}\star  (a\star  c)$$
$$ -2 (a\star  (b\star  c))\star  d + 2 (a\star  (b\star  d))\star  c
 +2 (b\star  (a\star  c))\star  d - 2 (b\star  (a\star  d))\star  c$$
$$ -2 (c\star  (a\star  d))\star  b + 2 (c\star  (b\star  d))\star  a
 +2 (d\star  (a\star  c))\star  b - 2 (d\star  (b\star  c))\star  a$$
$$ +2 ((a\star  c)\star  d)\star  b - 2 ((a\star  d)\star  c)\star  b
 -2 ((b\star  c)\star  d)\star  a + 2 ((b\star  d)\star  c)\star  a$$
 $$-((a\star  b)\star  c)\star  d + ((a\star  b)\star  d)\star  c - ((a\star  c)\star  b)\star  d+ 
  ((a\star  d)\star  b)\star  c$$
  $$
 +((b\star  a)\star  c)\star  d)- ((b\star  a)\star  d)\star  c+ ((b\star  c)\star  a)\star  d- 
  ((b\star  d)\star  a)\star  c$$
 $$+((c\star  a)\star  b)\star  d - ((c\star  b)\star  a)\star  d - ((d\star  a)\star  b)\star  c+ 
  ((d\star  b)\star  a)\star  c=0,$$
where we put $\{a,b\}=a\star b+b\star a.$
\end{itemize}

\section{$AR$-algebra of Zinbiel algebra}

\subsection{The multiplication $a\star b=\int a \int b$}
It is $AR$-multiplication for right-Zinbiel multiplication $a\star b=a\circ R(b)$ for $a\circ b=(\int a)b$ and $R(a)=\int b.$

{\bf Problem.} Is it true that any polynomial identity of algebra $(K[x],\star),$ where $a\star b=\int a\int b$ degree $k>2$ follows from commutativity identity? 

We have checked that until $k\le 6$ no identity except commutativity.

\subsection{The multiplication $a\star b=a\int_2 b$}

It is $AR$-multiplication for left-Zinbiel multiplication $a\star b=a\circ R(b)$ for $a\circ b=a\int b$ and $R(a)=\int b.$

The $AR$-multiplication satisfy the following identities. 
\begin{itemize}
    \item right-commutative,
    $$(a\star b)\star c=(a\star c)\star b$$
    \item  the following identity of degree $4$
    $$s_{1,3}(a,b,c,d)=(a\star [b, c])\star d+(a\star [c,d])\star b+(a\star [d, b])\star c,$$
    where $[a,b]=a\star b-b\star a.$
\item  the following identity of degree $4$
   $$ f_4'(a,b,c,d)=(a,b,[c,d])+(a,c,[d,b])+(a,d,[b,c])=0.$$
\end{itemize}

In particular, the algebra $(K[x],\star)^- =(K[x],[\;,\;])$, where $[a,b]=a\star b-b\star a=a \int_2 b-b\int_2 a$ is Tortkara. 
Under Jordan commutator $\{a,b\}=a\star b+b\star a=a \int_2 b+b\int_2 a$ the algebra $(K[x],\star)^+=(K[x],\{\;,\;\})$ does not staisfy the identity 
$f_5^+(a,b,c,d,e)=0$

\end{document}